\documentclass[12pt]{article}
\usepackage{amssymb}
\usepackage{amssymb,amsfonts}
\usepackage{amsmath}
\usepackage{cases}

\usepackage{indentfirst}
\usepackage[mathscr]{eucal}
\usepackage[OT1]{fontenc}

\usepackage{booktabs}

\usepackage[round,sort&compress]{natbib}
\usepackage{lscape}
\usepackage{multirow}
\usepackage{graphicx}
\usepackage{graphics}
\usepackage{color,xcolor}
\usepackage{float}

\usepackage{bm}

\usepackage[colorlinks,linkcolor=blue,anchorcolor=blue,citecolor=blue]{hyperref}

\usepackage{epstopdf}

\usepackage{tabls}
\usepackage{caption}
\usepackage{subcaption}
\usepackage{mathrsfs}
\usepackage{amsmath}
\allowdisplaybreaks[4]

\usepackage{geometry}
\newtheorem{theorem}{Theorem}

\newtheorem{lemma}{Lemma}[section]

\newtheorem{remark}{Remark}[section]

\usepackage[title]{appendix}

\numberwithin{equation}{section}

\def\Var{\mbox{Var}}
\def\Cov{\mbox{Cov}}
\def\E{\mbox{E}}
\def\tr{\mbox{tr}}

\def\C{{\bf {C}}}
\def\Z{{\bf {Z}}}
\def\B{{\bf {B}}}
\def\A{{\bf {A}}}
\def\M{{\bf {M}}}
\def\D{{\bf {D}}}
\def\Y{{\bf {Y}}}
\def\Q{{\bf {Q}}}
\def\S{{\bf {S}}}
\def\V{{\bf {V}}}
\def\I{{\bf {I}}}
\def\1{{\bf {1}}}
\def\0{{\bf {0}}}
\def\X{{\bf {X}}}

\def\P{{\bf {P}}}

\def\W{{\bf {W}}}
\def\e{{\bf {e}}}
\def\H{{\bf {H}}}

\def\G{{\bf {G}}}
\def\R{{\mathcal{R}}}

\newcommand{\bSigma}{\mbox{$\bm\Sigma$}}
\newcommand{\bmu}{\mbox{$\bm\mu$}}

\newcommand{\bGam}{\mbox{$\bm\Gamma$}}

\title{General linear hypothesis testing of high-dimensional mean vectors with unequal covariance matrices based on random integration}
\author{Mingxiang Cao$^{a}$,~Yelong Qiu$^{a}$ and Junyong Park$^{b*}$\\
    \it\footnotesize  $^{a}$Department of Statistics, Anhui Normal University, Wuhu, China\\
   \it\footnotesize   $^{b}$Department of Statistics, Seoul National University, Seoul, Korea}
\date{}

\begin{document}

\maketitle \footnotetext[1]{*Corresponding author.}

\footnotetext{E-mail:\url{junyongpark@snu.ac.kr}}

\noindent\textbf{Abstract}\ \ This paper is devoted to the study of the general linear hypothesis testing (GLHT) problem of multi-sample high-dimensional mean vectors. For the GLHT problem, we introduce a test statistic based on $L^2$-norm and random integration method, and deduce the asymptotic distribution of the statistic under given conditions. Finally, the potential advantages of our test statistics are verified by numerical simulation studies and examples.



\vspace{1cm}

\section{Introduction}\label{sec1}
\noindent  \quad With the rapid development of modern science and technology, we are faced with more and more high-dimensional data, which widely exist in many fields such as diffusion tensor imaging (\cite{Le Bihan et al.:2001}), DNA microarrays (\cite{Witten and Tibshirani:2011}), finance (\cite{Lam and Yao:2012}) and so on. Due to the explosion of dimensions, the study of statistical problems in the high-dimensional context will suffer great challenges, because the classical test methods are not applicable and often unstable in the high-dimensional framework. For example, in high-dimensional setting, since the sample covariance matrix is no longer invertible, the classical Hotelling's $T^2$-test is not applicable. In view of this situation, statisticians have done a lot of work in recent decades on high-dimensional hypothesis testing. Testing the equality of several mean vector for high-dimensional data is one of the important topics in hypothesis testing. In this paper, we are interested  testing problem on the homogeneity of multiple populations mean vectors in a high-dimensional context. Specifically, let $p$-dimensional random vectors $\Y_{i1},\ldots,\Y_{in_i}$ be independent sample from the $i$th population with the $p$-dimensional distribution $\Y_i$ with mean vector $\bmu_{i}$ and covariance matrix $\bSigma_{i}$ for $i=1,\ldots,K$, where dimension $p$ may be much larger than the overall sample size $n=\sum_{i=1}^{K}n_{i}$. we want to test the following general linear hypothesis testing problem without requiring $\bSigma_{i}$ to be equal for $i=1,\ldots,K$:
\begin{eqnarray}\label{H:1.1}
H_0:\widetilde{\G}\M=\0~~~~vs.~~~~~H_1:\widetilde{\G}\M\neq\0,
\end{eqnarray}
where $\M=(\bmu_{1},\bmu_{2},\cdots,\bmu_{K})^T$ is a $K\times p$ matrix and $\widetilde{\G}$ is  a $q\times K$ known coefficient matrix with rank$(\widetilde{\G})=q<K$. With different Settings of $\widetilde{\G}$, the hypothesis in (\ref{H:1.1}) can generate some special hypothesis testing problems. For example, when $\widetilde{\G}= (\I_{k-1},-\1_{k-1})$ or $\widetilde{\G}= (-\1_{k-1},\I_{k-1})$, the hypothesis in \eqref{H:1.1} can generate the following hypothesis testing problem:
\begin{eqnarray}\label{H:1.2}
H_0:\bmu_{1}=\bmu_{2}=\cdots=\bmu_{K}~~~~vs.~~~~~H_1:H_0\mbox{ is not true},
\end{eqnarray}
where $\I_{K-1}$ and $\1_{K-1}$ denote the $(K-1)\times (K-1)$ identity matrix and the $(K-1)$-dimensional vector with all elements $1$, respectively. the hypothesis in (\ref{H:1.2}) is the well-known multivariate analysis of variance (MANOVA) problem, which has been extensively studied in a large literature. For example, \cite{Cai and Xia:2014} and \cite{Chakraborty and Sakhanenko:2023}. If we take $\widetilde{\G}$ as a $K$-dimensional row vector ($c_{1},\ldots,c_{K}$), then the hypothesis in \eqref{H:1.1} contains the following hypothesis testing problem about linear combination of $K$ mean vectors:
\begin{eqnarray}\label{H:1.3}
H_0:\sum_{i=1}^{K}{c_{i}\bmu_{i}}=\0~~~~vs.~~~~~H_1:\sum_{i=1}^{K}{c_{i}\bmu_{i}}\neq\0,
\end{eqnarray}
this particular test problem is studied in \cite{Nishiyama et al.:2013} and \cite{Cao et al.:2022} respectively.

The hypothesis in \eqref{H:1.1} is equivalently rewritten as 
\begin{eqnarray}\label{H:1.4}
H_0:\G\M=\0~~~~vs.~~~~~H_1:\G\M\neq\0,
\end{eqnarray}
where $\G=(\widetilde{\G}\D\widetilde{\G}^T)^{1/2}\widetilde{\G}$ and $\D=diag(\frac{n}{n_1},\ldots,\frac{n}{n_K})$. The advantage of the above transformation is that the test statistic proposed for \eqref{H:1.1} is invariant under this non-singular transformation $\widetilde{\G}\rightarrow\P\widetilde{\G}$, where $\P$ is an arbitrary $q\times q$ non-singular matrix.

\begin{remark}
This matrix $\D$ is related to the covariance matrix of $\sqrt{n}\widehat{\bmu}$, where $\widehat{\bmu}=(\overline{\Y}_{1}^T,\ldots,\overline{\Y}_{K}^T)^T$ is an unbiased estimate of $\bmu=(\bmu_1^T,\ldots,\bmu_K^T)^T$. We can get $Cov(\sqrt{n}\widehat{\bmu})=diag(\frac{n}{n_1}\bSigma_{1},\ldots,\frac{n}{n_K}\bSigma_{K})$.
\end{remark}

Further, we can write \eqref{H:1.4} as the following equivalent form:
\begin{eqnarray}\label{H:1.5}
H_0:\C\bmu=\0~~~~vs.~~~~~H_1:\C\bmu\neq\0,
\end{eqnarray}
where $\C=\G\otimes\I_p$ and $\otimes$ represents the Kronecker product operator. This hypothesized form has been used in many literature, such as \cite{Zhou et al.:2017}, \cite{Zhang et al.:2022}, \cite{Zhang and Zhu:2022a} and so on.

When $K=2$, the hypothesis in \eqref{H:1.2} will be reduced to a two-sample problem in high-dimensional context. Studies of the two-sample problem can be traced back to \cite{Dempster:1958}, which proposed a non-exact tests through the $F$-distribution approximation. Subsequently, \cite{Bai and Saranadasa:1996} proposed a non-exact test by modifying Hotelling's $T^2$ test. However their tests assume that the covariance matrix of the two samples is the same, but when the two covariance matrix are not equal, it will produce misleading results. To overcome this difficulty, \cite{Chen and Qin:2010} proposed a test based on U-statistic when the sample covariance matrix was not equal. Since then, a large number of different tests have been proposed, such as \cite{Srivastava et al.:2013}, \cite{Feng and Sun:2015}, \cite{Wang and Xu:2022} and so on.

Furthermore, the research of these pioneers has been extended to high-dimensional $K$-samples, and there are many research works available. For $K>2$, when $K$ samples have equal covariance matrices, \cite{Schott:2007} extended the work of \cite{Bai and Saranadasa:1996} to a MANOVA test, \cite{Srivastava and Kubokawa:2013} proposed a scale-invariant MANOVA test for non-normal data, and \cite{Jiang et al.:2015} proposed a non parametric $K$-sample test. When $K$ samples have different covariance matrices, \cite{Yamada and Himeno:2015} proposed a MANOVA test based on the variation matrix, \cite{Hu et al.:2017} were inspired by the thought of \cite{Chen and Qin:2010} and extended the two sample test to the MANOVA test, and \cite{Li et al.:2024} proposed a test for linear combinations of $K$ mean vectors based on the $L^2$-norm. \cite{Zhou et al.:2017} proposed a test for the GLHT problem based on the U-statistic, and \cite{Zhang et al.:2022} tested the GLHT problem based on the $L^2$-norm. 

In this article, using the random integration technique proposed by \cite{Jiang et al.:2024}, we propose a test for GLHT problem based on U-statistic, without assuming the normality of the $K$ population or the same $K$ covariance matrices. Although our test statistic is based on the U-statistic, we can quickly calculate it without using the knowledge of the U-statistic. We establish the asymptotic distribution under mild conditions, and give the sample ratio-consistent estimator. In addition, we demonstrate the validity of our proposed statistics through simulation studies.

The remainder of this article is organized as follows. In Section \ref{sec2}, we first derive the test statistic using the random integration method and then present our main results. In Section \ref{sec3}, a simulation study is presented. Section \ref{sec4} presents a real data simulation. In Section \ref{sec5} provides the conclusion. All technical proofs are arranged in the Appendix.

\section{Test procedure and main results}\label{sec2}
\subsection{Test procedure}\label{sec2.1}
\noindent \quad To propose our test statistic, we apply the random integration method to the GLHT problem. Therefore, to test whether the null hypothesis $\C\bmu=0$ in \eqref{H:1.5} holds, we assume $\Y=(\Y_1^T,\ldots,\Y_K^T)^T$ and derive the following equivalence
\begin{align}\label{eq2.1} 
 \C\bmu=0&\Leftrightarrow\E(\C\Y)=0\nonumber\\
&\Leftrightarrow\E[(\I_q\otimes\delta^T)\C\Y]=\E[(\G\otimes\delta^T)\Y]=0 \quad\forall\delta\in\mathbb{\R}^p \nonumber \\
&\Leftrightarrow\E^{2}[(\G\otimes\delta^T)\Y]=0\quad\forall\delta\in\mathbb{\R}^p,
\end{align}
where $\Y_i$ is the distribution of the $i$-th population with mean vector $\bmu_{i}$ and covariance matrix $\bSigma_{i}$ for $i=1,\ldots,K$.  
For convenience, set $\G=(\bf g_1,\ldots,\bf g_{K})$, $\D=(\G\otimes\delta^T),$ so we have 
\begin{align}\label{eq2.2} 
\D=(\D_1,\ldots,\D_K), \quad\D_{\alpha}^T\D_{\beta}=d_{\alpha\beta}\delta\delta^T, \quad\G^T\G=\left(d_{\alpha\,\beta}\right)_{\alpha,\beta=1}^K,
 \end{align}
where $\D_{\alpha}={\bf g_{\alpha}}\otimes\delta^{T}$, $d_{\alpha\beta}={\bf g_\alpha}^T{\bf g_\beta}$, $\alpha,\beta=1,\ldots,K.$
It shows from these equation in \eqref{eq2.1} that testing whether the null hypothesis $H_0$ in \eqref{H:1.5} holds amounts to testing whether
\begin{align}\label{eq2.3} 
{\rm{RI}}_{w}:=\int{\E^{2}[(\G\otimes\delta^T)\Y]}\omega(\delta)d\delta=\0,
 \end{align}
where $\omega(\delta)$ denotes any positive weight function. To obtain the explicit expression of \eqref{eq2.3}, we derive the following theorem \ref{th2.1}.

\begin{theorem}\label{th2.1}
Let $\omega(\delta)=\prod_{i=1}^{p}\omega_i(\delta_i)$, where $\omega_i(\cdot)$ is a density function with mean $\alpha_i$ and variance $\beta_i^2$ for $i=1,\ldots,p$, then
\begin{align}\label{eq2.4}
{\rm{RI}}:={\rm{RI}}_{w}=\sum_{\alpha=1}^{K}d_{\alpha\alpha}\bmu_\alpha^T\W\bmu_\alpha+\sum_{\alpha\neq\beta}d_{\alpha\beta}\bmu_\alpha^T\W\bmu_\beta,
\end{align}
and ${\rm{RI}}\geq0$ with the equality holds if and only if  $\C\bmu=\0$, where $\W=\B+{\bf a}{\bf a}^T, {\bf a}=(\alpha_1,\ldots,\alpha_p)^T, \B=diag(\beta_1^2,\ldots,\beta_p^2)$.
\end{theorem}

The expression in Theorem \ref{th2.1} indicates that the null hypothesis in \eqref{H:1.5} becomes the equivalent expression $\rm{RI}=0$. Therefore, based on the explicit form of Theorem \ref{th2.1}, we can easily see that for $\alpha,\beta=1,2,\ldots,K$, the U-statistic for estimating $\bmu_{\alpha}^T\W\bmu_\alpha$ and $\bmu_{\alpha}^T\W\bmu_\beta$ can be respectively given by 
\begin{align}\label{eq2.5}
\S_{\alpha\alpha}=\frac{\sum\limits_{i\neq j}\Y_{\alpha i}^T\W\Y_{\alpha j}}{n_{\alpha}(n_{\alpha}-1)} \quad and\quad \S_{\alpha\beta}=\frac{\sum\limits_{i,j}\Y_{\alpha i}^T\W\Y_{\beta j}}{n_{\alpha}n_{\beta}}.
\end{align}
Then for the test ${\rm{RI}}$, we propose a new test statistic based on the U-statistic:
\begin{align}\label{eq2.6}
{\rm{T}}_n=\sum_{\alpha=1}^{K}d_{\alpha\alpha}\S_{\alpha\alpha}+\sum_{\alpha\neq\beta}d_{\alpha\beta}\S_{\alpha\beta},
\end{align}
which has expectation ${\rm{RI}}$. 
\begin{remark}
When we test the hypothesis problem of \eqref{H:1.3}, our test statistic are similar to that of \cite{Li et al.:2024}, but the problem we study is more general.
\end{remark}

By simple linear algebra calculations, for $\alpha,\beta=1,2,\ldots,K$, we have
\begin{align}\label{eq2.7}
\S_{\alpha\alpha}=\overline{\Y}_{\alpha}^T\W\overline{\Y}_{\alpha}-\frac{{\rm{tr}}(\W\widehat{\bSigma_{\alpha}})}{n_\alpha},
\quad \S_{\alpha\beta}=\overline{\Y}_{\alpha}^T\W\overline{\Y}_{\beta},
\end{align}
where 
\begin{align}\label{eq2.8}
\overline{\Y}_{\alpha}=\frac{1}{n_\alpha}\sum\limits_{i=1}^{n_\alpha}\Y_{\alpha i},
\quad \widehat{\bSigma_{\alpha}}=\frac{1}{n_\alpha-1}\sum\limits_{i=1}^{n_\alpha}(\Y_{\alpha i}-\overline{\Y}_{\alpha})(\Y_{\alpha i}-\overline{\Y}_{\alpha})^T
\end{align}
are the sample mean and covariance matrix of the $\alpha$-th sample, respectively. After simple calculations from \eqref{eq2.5}, \eqref{eq2.6} and \eqref{eq2.7}, we can obtain 
\begin{align}\label{eq2.9}
{\rm{T}}_n&=\sum_{\alpha,\beta}d_{\alpha\beta}\overline{\Y}_{\alpha}^T\W\overline{\Y}_{\beta}-
\sum_{\alpha=1}^{K}d_{\alpha\alpha}\frac{{\rm{tr}}(\W\widehat{\bSigma_{\alpha}})}{n_\alpha}\\
&={\lVert(\G\otimes\W^{\frac{1}{2}})\widehat{\bmu }\rVert}^2-\sum_{\alpha=1}^{K}d_{\alpha\alpha}\frac{{\rm{tr}}(\W\widehat{\bSigma_{\alpha}})}{n_\alpha},
\end{align}
where $\widehat{\bmu}=(\overline{\Y}_{1}^T,\ldots,\overline{\Y}_{K}^T)^T$ is an unbiased estimate of $\bmu=(\bmu_1^T,\ldots,\bmu_K^T)^T$. The above equation indicates that we can quickly calculate our test statistic ${\rm{T}}_n$ without utilizing knowledge of the U-statistic. Further, let $\D_{\theta}=(\G\otimes\W^{\frac{1}{2}})$, we can write 
\begin{align}\label{eq2.10}
{\rm{T}}_n={\rm{T}}_{n_0}+2\S_{n}+{\lVert\D_{\theta}\bmu\rVert}^2,
\end{align}
where 
\begin{align}\label{eq2.11}
{\rm{T}}_{n_0}={\lVert\D_{\theta}(\widehat{\bmu}-\bmu)\rVert}^2-\sum_{\alpha=1}^{K}\frac{d_{\alpha\alpha}}{{n_\alpha}}{\rm{tr}}(\W\widehat{\bSigma_{\alpha}}),
\quad\S_n=(\D_{\theta}\bmu)^{T}\D_{\theta}(\widehat{\bmu}-\bmu).
\end{align}
For the convenience of further discussion, let
\begin{align}\label{eq2.12}
\X_{\alpha i}=\Y_{\alpha i}-\bmu_{\alpha},\quad i=1,\ldots,n_{\alpha};\quad \alpha=1,\ldots,K.
\end{align}
Therefore, the $K$ samples mentioned above are independent, and for each $\alpha=1,\ldots,K$, $\X_{\alpha 1},\ldots,\X_{\alpha n_{\alpha}}$ is i.i.d. with $\E(\X_{\alpha 1})=0$ and $\Cov(\X_{\alpha 1})=\bSigma_{\alpha}$. Furthermore, we have
\begin{align*}
\overline{\X}_{\alpha} &=\frac{1}{n_\alpha}\sum_{i=1}^{n_\alpha}\X_{\alpha i}=\overline{\Y}_{\alpha}-\bmu_{\alpha},\\
\widehat{\bSigma_{\alpha}}&=\frac{1}{n_\alpha-1}\sum_{i=1}^{n_\alpha}(\X_{\alpha i}-\overline{\X}_{\alpha})(\X_{\alpha i}-\overline{\X}_{\alpha})^T\\ 
&= \frac{1}{n_\alpha-1}\sum_{i=1}^{n_\alpha}(\Y_{\alpha i}-\overline{\Y}_{\alpha})(\Y_{\alpha i}-\overline{\Y}_{\alpha})^T.                                            
\end{align*}
From these expressions, together with \eqref{eq2.6} and \eqref{eq2.9}, we can derive
\begin{align}\label{eq2.13}
{\rm{T}}_{n_0}=\sum_{\alpha,\beta}d_{\alpha\beta}\overline{\X}_{\alpha}^T\W\overline{\X}_{\beta}-
\sum_{\alpha=1}^{K}d_{\alpha\alpha}\frac{{\rm{tr}}(\W\widehat{\bSigma_{\alpha}})}{n_\alpha}
=\sum_{\alpha=1}^{K}d_{\alpha\alpha}\S_{\alpha\alpha}^{0}+\sum_{\alpha\neq\beta}d_{\alpha\beta}\S_{\alpha\beta}^{0},
\end{align}
where 
\begin{align}\label{eq2.14}
\S_{\alpha\alpha}^{0}=\frac{\sum\limits_{i\neq j}\X_{\alpha i}^T\W\X_{\alpha j}}{n_{\alpha}(n_{\alpha}-1)},\quad \S_{\alpha\beta}^{0}=\frac{\sum\limits_{i,j}\X_{\alpha i}^T\W\X_{\beta j}}{n_{\alpha}n_{\beta}}.
\end{align}

Under $\H_{0}$, $\S_{n}=0$ and ${\lVert\D_{\theta}\bmu\rVert}^2$ is equal to 0. So ${\rm{T}}_{n_0}$ has the same distribution as ${\rm{T}}_{n}$ under the null hypothesis. Furthermore, we can use the expressions in \eqref{eq2.13} and \eqref{eq2.14} to study the null distribution of ${\rm{T}}_{n}$. Let $\alpha$, $\beta$ and $\gamma$ be different integers, then by simple calculation we can get 
\begin{align}\label{eq2.15}
&\E(\S_{\alpha\alpha}^{0})=\E(\S_{\alpha\beta}^{0})=0,\nonumber\\
&\Var(\S_{\alpha\alpha}^{0})=\frac{2{\rm{tr}}(\W\bSigma_{\alpha})^{2}}{n_{\alpha}(n_{\alpha}-1)},\quad
\Var(\S_{\alpha\beta}^{0})=\frac{{\rm{tr}}(\W\bSigma_{\alpha}\W\bSigma_{\beta})}{n_{\alpha}n_{\beta}},\\
&\Cov(\S_{\alpha\alpha}^{0},\S_{\beta\beta}^{0})=\Cov(\S_{\alpha\beta}^{0},\S_{\beta\beta}^{0})=\Cov(\S_{\alpha\beta}^{0},\S_{\alpha\gamma}^{0})=0\nonumber.
\end{align}
Therefore, under the null hypothesis, the variance of T can be given by 
\begin{align}\label{eq2.16}
\Var({\rm{T}}_{n_0})=2\left(\sum_{\alpha=1}^{K}\frac{d_{\alpha\alpha}^{2}{\rm{tr}}(\W\bSigma_{\alpha})^{2}}{n_{\alpha}(n_{\alpha}-1)}+
\sum_{\alpha\neq\beta}\frac{d_{\alpha\beta}^{2}{\rm{tr}}(\W\bSigma_{\alpha}\W\bSigma_{\beta})}{n_{\alpha}n_{\beta}}\right). 
\end{align}

\begin{remark}
The expression of ${\rm{T}}_{n}$ in \eqref{eq2.9} can quickly calculate ${\rm{T}}_{n}$ without utilizing knowledge of the U-statistic. However, the expression of ${\rm{T}}_{n_0}$ based on the U-statistic in \eqref{eq2.13} can quickly derive the variance of ${\rm{T}}_{n}$ under $\H_{0}$. This makes our calculations more efficient.
\end{remark}

\subsection{Main results}\label{sec2.2}
\noindent \quad In this section, we will investigate the asymptotic normality of our test ${\rm{T}}_{n}$. For this purpose, we impose the following assumptions:
\begin{enumerate}
\item[(A1)] We assume that the samples are generated from the multivariate factor model $\Y_{\alpha i}=\bGam_{\alpha}\Z_{\alpha i}+\bmu_{\alpha}$ for $i=1,\ldots,n_{\alpha}$ and $\alpha=1,\ldots,K$, where $\bGam_{\alpha}$ is a $p\times m$ matrix with $p\leq m$ satisfying $\bGam_{\alpha}\bGam_{\alpha}^T=\bSigma_{\alpha}$, and $\Z_{\alpha i}$ are $i.i.d.$ $m$-dimensional vector with $\E(\Z_{\alpha i})=0$ and $\Cov(\Z_{\alpha i})=\I_m$. Then \eqref{eq2.12} can be rewritten as $\X_{\alpha i}=\bGam_{\alpha}\Z_{\alpha i}$, for $i=1,\ldots,n_{\alpha}$ and $\alpha=1,\ldots,K$.
\item[(A2)] Let $\E(z_{\alpha i s}^4)=3+\Delta<\infty$, where $z_{\alpha i s}$ is the $s$-th component of $\Z_{\alpha i}$ for $i=1,\ldots,n_{\alpha}$, $\alpha=1,\ldots,K$, and $\Delta$ is some constant. Further, 
    \begin{align}\label{eq2.17}
\E(z_{\alpha is_{1}}^{\varsigma_{1}}\cdots{z}_{\alpha is_{r}}^{\varsigma_{r}})=\E(z_{\alpha is_{1}}^{\varsigma_{1}})\cdots\E(z_{\alpha is_{r}}^{\varsigma_{r}})
\end{align}
for any positive integers $r$ and $\varsigma_{s}$ where $\sum_{s=1}^{r}\varsigma_{s}\leq 8$, $s_{1},\ldots,s_{r}$ are distinct indices.
\item[(A3)] $n_{\alpha}/n\to\tau_{\alpha}\in(0,1)$ as $n\to\infty$, where $\alpha=1,\ldots,K$, $n=n_{1}+\cdots+n_{K}$.
\item[(A4)] As $p\to\infty$, for any $\alpha,\beta,\gamma\in\{1,\ldots,K\}$,
\begin{align}\label{eq2.18}
\tr(\W\bSigma_{\alpha}\W\bSigma_{\beta}\W\bSigma_{\gamma}\W\bSigma_{\beta})=
o\left\{\tr(\W\bSigma_{\alpha}\W\bSigma_{\beta})\tr(\W\bSigma_{\gamma}\W\bSigma_{\beta})\right\}.
\end{align}

\end{enumerate}

\begin{remark}
The above assumptions are derived by extending the conditions in \cite{Chen and Qin:2010} and \cite{Jiang et al.:2024} to the background of our article. It shows that we do not require a clear relationship between the data dimension $p$ and the total sample size $n$. Assumption (A1) and (A2) are widely used multivariate factor models in high-dimensional data analysis literature. The formula in \eqref{eq2.17} requires that the the components of $\Z_{\alpha i}$ have a pseudo-independence. Assumption (A3) is the standard regularity hypothesis in high-dimensional multi-sample tests, which ensures that the group sample sizes tend proportionally to infinity. Assumption (A4) ensures the consistency and asymptotic normality of our test.

\end{remark}

To obtain the asymptotic normality of ${\rm{T}}_{n_0}$, we need to construct square integrable martingale with zero mean. For convenience, we define $\Y_{m_{\alpha-1}+i}=\W^{1/2}\X_{\alpha i}$, where $m_0=0,m_{\alpha}=n_1+\ldots+n_{\alpha}$, $i=1,\ldots,n_{\alpha}$ and $\alpha=1,\ldots,K$. let
\begin{align}\label{eq2.19}
\phi_{ij}= \begin{cases}
\frac{2d_{\alpha\beta}}{n_{\alpha}n_{\beta}}\Y_{i}^{T}\Y_{j}, & i\in \{m_{\alpha-1}+1,\ldots,m_{\alpha}\},j\in \{m_{\beta-1}+1,\ldots,m_{\beta}\},\alpha<\beta, \\
\frac{2d_{\alpha\alpha}}{n_{\alpha}(n_{\alpha}-1)}\Y_{i}^{T}\Y_{j}, & i,j\in \{m_{\alpha-1}+1,\ldots,m_{\alpha}\},i<j.
\end{cases}
\end{align}
Then ${\rm{T}}_{n_0}$ in \eqref{eq2.13} can be written as 
\begin{align}\label{eq2.20}
{\rm{T}}_{n_0}&=\sum_{\alpha=1}^{K}\frac{2d_{\alpha\alpha}}{n_{\alpha}(n_{\alpha}-1)}\sum_{i<j}\X_{\alpha i}^T\W\X_{\alpha j}+
\sum_{1\leq\alpha<\beta\leq K}\frac{2d_{\alpha\beta}}{n_{\alpha}n_{\beta}}\sum_{i,j}\X_{\alpha i}^T\W\X_{\beta j}\nonumber\\
&=\sum_{j=2}^{m_K}\sum_{i=1}^{j-1}\phi_{ij}.
\end{align}
Let $\V_{nj}=\sum_{i=1}^{j-1}\phi_{ij}$ and $\S_{nm}=\sum_{j=2}^{m}\sum_{i=1}^{j-1}\phi_{ij}=\sum_{j=2}^{m}\V_{nj}$, $m=2,\ldots,m_K$ and $\mathcal{F}_{nm}=\sigma\{\Y_1,\ldots,\Y_m\}$, which is the $\sigma$-ﬁeld generated by $\{\Y_1,\Y_2,\ldots,\Y_m\}$. Then, we have the following lemma.
\begin{lemma}\label{lem 2.1}
For each $n$, $\{\S_{nm}\}_{m=1}^{m_K}$ is a zero mean, square integrable martingale about the filtration $\{\mathcal{F}_{nm}\}_{m=1}^{m_K}$.
\end{lemma}

Let $\V_{n1}=0$, $\mathcal{F}_{n0}=\{\emptyset,\Omega\}$, and $\sigma({\rm{T}}_{n_0})=\sqrt{\Var({\rm{T}}_{n_0})}$. To facilitate the proof of asymptotic normality of ${\rm{T}}_{n_0}$, we derive the following lemmas.
\begin{lemma}\label{lem 2.2}
Under Assumption (A1)-(A4), for $i\in\{m_{\gamma-1}+1,\ldots,m_{\gamma}\}$, we have
\begin{align}\label{eq2.21}
&\E(\Y_{i}^{T}\W^{1/2}\bSigma_{\alpha}\W^{1/2}\Y_{i}\Y_{i}^{T}\W^{1/2}\bSigma_{\beta}\W^{1/2}\Y_{i})\nonumber\\
&\leq \tr(\W\bSigma_{\alpha}\W\bSigma_{\gamma})\tr(\W\bSigma_{\beta}\W\bSigma_{\gamma})
+(2+\Delta)\tr(\W\bSigma_{\alpha}\W\bSigma_{\gamma}\W\bSigma_{\beta}\W\bSigma_{\gamma}).
\end{align}
\end{lemma}
Furthermore, for $i\in\{m_{\alpha-1}+1,\ldots,m_{\alpha}\}$, $j\in\{m_{\beta-1}+1,\ldots,m_{\beta}\}$, $k\in\{m_{\gamma-1}+1,\ldots,m_{\gamma}\}$, where $i\neq j$, $k\neq i$, and $k\neq j$, we have $\E\left((\Y_{i}^{T}\Y_{k})^{2}(\Y_{j}^{T}\Y_{k})^{2}\right)$ not greater than the right-hand side of (2.22).

\begin{lemma}\label{lem 2.3}
Under Assumption (A1)-(A4), for $i\in\{m_{\alpha-1}+1,\ldots,m_{\alpha}\}$, $j\in\{m_{\beta-1}+1,\ldots,m_{\beta}\}$, and $i\neq j$, we have
\begin{align}\label{eq2.22}
\E(\Y_{i}^{T}\Y_{j})^{4}\leq (3+\Delta)\tr^{2}(\W\bSigma_{\alpha}\W\bSigma_{\beta})+(3+\Delta)(2+\Delta)\tr(\W\bSigma_{\alpha}\W\bSigma_{\beta})^{2}.
\end{align}
\end{lemma}

\begin{lemma}\label{lem 2.4}
Under Assumption (A1)-(A4), we have
\begin{align}\label{eq2.23}
\Var\left(\sum_{j=m_{\alpha-1}+1}^{m_{\alpha}}\E(\V_{nj}^{2}|\mathcal{F}_{n,j-1})\right)=o(\sigma^{4}({\rm{T}}_{n_0})), \quad\alpha=1,\ldots,K,\nonumber\\
\sum_{j=1}^{n}\E(\V_{nj}^{2}|\mathcal{F}_{n,j-1})/\sigma^{2}({\rm{T}}_{n_0})\stackrel{p}{\longrightarrow}1,\quad \sum_{j=1}^{n}\E(\V_{nj}^{4})=o(\sigma^{4}({\rm{T}}_{n_0})).
\end{align}
\end{lemma}

Based on these lemmas and the martingale difference central limit theorem, we can establish the asymptote normality of ${\rm{T}}_{n_0}$ in the following theorem.
\begin{theorem}\label{th2.2}
Under Assumption (A1)-(A4), as $n,p\to\infty$,
$$\frac{{\rm{T}}_{n_0}}{\sigma({\rm{T}}_{n_0})}\stackrel{d}{\longrightarrow}N(0,1).$$
\end{theorem}

In order to construct a test procedure, we need to derive the ratio-consistent estimator of $\sigma^{2}({\rm{T}}_{n_0})$. Obviously the unbiased estimator  of ${\rm{tr}}(\W\bSigma_{\alpha}\W\bSigma_{\beta})$ can be given by ${\rm{tr}}(\W\widehat{\bSigma_{\alpha}}\W\widehat{\bSigma_{\beta}})$. Since the estimator for ${\rm{tr}}(\W\bSigma_{\alpha})^{2}$ proposed by \cite{Chen and Qin:2010} is biased as indicated by \cite{Feng et al.:2015}, We use the unbiased estimator method proposed by \cite{Himeno and Yamada:2014} for ${\rm{tr}}(\W\bSigma_{\alpha})^{2}$ to derive the unbiased estimator of ${\rm{tr}}(\W\bSigma_{\alpha})^{2}$. Using the idea of replacing $\Y$ with $\W^{1/2}\Y$, we provide an unbiased estimator of ${\rm{tr}}(\W\bSigma_{\alpha})^{2}$ by
\begin{align*}
{\rm{tr}}\widehat{(\W\bSigma_{\alpha})^2}=\frac{n_{\alpha}-1}{n_{\alpha}(n_{\alpha}-2)(n_{\alpha}-3)}\left((n_{\alpha}-1)(n_{\alpha}-2){\rm{tr}}(\W\widehat\bSigma_{\alpha})^{2
}+{\rm{tr}}^{2}(\W\widehat\bSigma_{\alpha})-n_{\alpha}\Q_{\alpha}\right),
\end{align*}
where $\Q_{\alpha}=(n_{\alpha}-1)^{-1}\sum_{i=1}^{n_{\alpha}}{\lVert\W^{1/2}(\Y_{\alpha i}-\overline{\Y}_{\alpha})\rVert}^4$. Thus, from \eqref{eq2.16}, the unbiased estimator of $\sigma^{2}({\rm{T}}_{n_0})$ can be given by 
\begin{align}\label{eq2.24}
\widehat\sigma^{2}({\rm{T}}_{n_0})=2\left(\sum_{\alpha=1}^{K}\frac{d_{\alpha\alpha}^{2}{\rm{tr}}\widehat{(\W\bSigma_{\alpha})^{2}}}{n_{\alpha}(n_{\alpha}-1)}+
\sum_{\alpha\neq\beta}\frac{d_{\alpha\beta}^{2}{\rm{tr}}(\W\widehat{\bSigma_{\alpha}}\W\widehat{\bSigma_{\beta}})}{n_{\alpha}n_{\beta}}\right).
\end{align}
The advantage of equation \eqref{eq2.24} is that we can quickly calculate the variance estimator of our test statistic without utilizing knowledge of the U-statistic. Further, we prove that our proposed unbiased estimators are ratio-consistent, as shown in the following theorem.

\begin{theorem}\label{th2.3}
Under Assumption (A1)-(A4), as $n,p\to\infty$, for $\alpha,\beta=1,\ldots,K, \alpha\neq\beta$, we have
$$\frac{{\rm{tr}}\widehat{(\W\bSigma_{\alpha})^2}}{{\rm{tr}}(\W\bSigma_{\alpha})^{2}}\stackrel{p}{\longrightarrow}1,\quad
\frac{{\rm{tr}}(\W\widehat{\bSigma_{\alpha}}\W\widehat{\bSigma_{\beta}})}{{\rm{tr}}(\W\bSigma_{\alpha}\W\bSigma_{\beta})}\stackrel{p}{\longrightarrow}1,\quad
\frac{\widehat\sigma^{2}({\rm{T}}_{n_0})}{\sigma^{2}({\rm{T}}_{n_0})}\stackrel{p}{\longrightarrow}1.$$
\end{theorem}

Furthermore, by Theorem \ref{th2.2}, Theorem \ref{th2.3}, and Slutsky’s theorem, under Assumption (A1)-(A4) and $H_0$, as $n,p\to\infty$, we have $\frac{{\rm{T}}_{n}}{\widehat\sigma({\rm{T}}_{n_0})}\stackrel{p}{\longrightarrow}N(0,1)$, where $\widehat\sigma({\rm{T}}_{n_0})=\sqrt{\widehat\sigma^{2}({\rm{T}}_{n_0})}$.

Next, we will investigate the asymptotic power of ${\rm{T}}_{n}$. From \eqref{eq2.10}, we have ${\rm{T}}_n={\rm{T}}_{n_0}+2\S_{n}+{\lVert\D_{\theta}\bmu\rVert}^2$, where ${\rm{T}}_{n_0}$ has the same distribution as ${\rm{T}}_{n}$ under $H_{0}$, and $\S_n=(\D_{\theta}\bmu)^{T}\D_{\theta}(\hat{\bmu}-\bmu)$. By calculation, we have 
\begin{align}\label{eq2.25}
\Var(\S_n)&=(\D_{\theta}\bmu)^{T}\D_{\theta}diag\left(\frac{\bSigma_{1}}{n_1},\ldots,\frac{\bSigma_{K}}{{n_K}}\right)\D_{\theta}^{T}(\D_{\theta}\bmu)\nonumber\\
&=\sum_{\alpha=1}^{K}\sum_{\beta=1}^{K}\sum_{\gamma=1}^{K}\frac{d_{\alpha\beta}d_{\beta\gamma}\bmu_{\alpha}^{T}\W\bSigma_{\beta}\W\bmu_{\gamma}}{n_{\beta}}.
\end{align}
Furthermore, to obtain the asymptotic power of ${\rm{T}}_{n}$ under the non-null hypothesis, we use the following condition 
\begin{align}\label{eq2.26}
\Var(\S_n)=o(\sigma^{2}({\rm{T}}_{n_0})),
\end{align}
which describes the situation where the information in local alternative is relatively small compared to the variance of T under the null hypothesis. That is, $|\Var({\rm{T}}_{n})-\Var({\rm{T}}_{n_0}|\rightarrow 0$ as $n\rightarrow\infty$. Therefore, we can easily obtain the asymptotic power of ${\rm{T}}_{n}$, which is given in the following theorem.
\begin{theorem}\label{th2.4}
Under Assumption (A1)-(A4) and \eqref{eq2.16} as $n,p\to\infty$, we have
\begin{align}\label{eq2.27}
P\left(\frac{{\rm{T}}_{n}}{\widehat\sigma({\rm{T}}_{n_0})}\geq z_{\alpha}\right)=\Phi\left(-z_{\alpha}+\frac{{\lVert\D_{\theta}\bmu\rVert}^2}{\sigma({\rm{T}}_{n_0})}\right)(1+o(1)) ,
\end{align}
\end{theorem}
where $z_{\alpha}$ and $\Phi(\cdot)$ denote the upper $100\alpha$ percentile and the cumulative distribution function of the standard normal distribution, respectively.

\begin{remark}
The the asymptotic power of ${\rm{T}}_{n}$ proposed by us is mainly determined by $\frac{{\lVert\D_{\theta}\bmu\rVert}^2}{\sigma({\rm{T}}_{n_0})}$. If this ratio is positive, the asymptotic power will be extraordinary, and if this ratio tends towards $\infty$, the asymptotic power will tend towards $1$.
\end{remark}

\section{Simulation study}\label{sec3}
\noindent \quad In this section, we will evaluate the performance of our proposed test $T_{New}$ through two simulation studies. Simulation 1 and Simulation 2 respectively aim to compare the performance of our test with three existing tests for one-way MANOVA and GLHT problems under different configurations. For convenience, we use $T_{Z}$, $T_{ZZG}$, and $T_{ZZ}$ to represent the tests of \cite{Zhou et al.:2017}, \cite{Zhang et al.:2022}, and \cite{Zhang and Zhu:2022a}, respectively.

Further, in all simulations, these $K$ high-dimensional samples are generated using the factor model in Assumption (A1): 
$$\Y_{\alpha i}=\bGam_{\alpha}\Z_{\alpha i}+\bmu_{\alpha}\quad for\quad i=1,\ldots,n_{\alpha} \quad and \quad\alpha=1,\ldots,K ,$$
where $\Z_{\alpha i}=(z_{\alpha i 1},\ldots,z_{\alpha i p})^T$ and $z_{\alpha i k}, k=1,2,\dots,p$ is i.i.d. and generated from the following three distributions: 
\begin{itemize}
\item[] Model 1: standard normal distribution $\mathcal{N}(0,1)$.

\item[] Model 2: standardized $t$ distribution $t_{4}/\sqrt{2}$.

\item[] Model 3: standardized chi-square distribution $(\chi_{1}^{2}-1)/\sqrt{2}$.
\end{itemize}

The empirical size and power of all simulations are calculated from $2000$ replicates performed at the normal level $\alpha=0.05$. In the simulation, we set $K=4$, sample size ${\bf n}=(n_{1},n_{2},n_{3},n_{4})=(20,30,45,50), (35,60,80,90), (40,70,100,120)$ and data dimension $p=100, 200, 300, 400, 500$. For the tuning parameter of our test, we take $\alpha_{1}=\alpha_{2}=\cdots=\alpha_{p}=2p^{-3/8}$ and $\beta_{i}=\sqrt{2}(p+i)/p, i=1,2,\ldots,p$. For the null hypothesis, we assume $\bmu_{1}=\bmu_2=\bmu_3=\bmu_4=0$. For the alternative hypothesis, we choose $\bmu_{1}=\bmu_2=\bmu_3=0$ and $\bmu_4$ has $[p^{1-t}]$ non-zero terms of equal value $\sqrt{2r(1/n_{1}+1/n_{2}+1/n_{3}+1/n_{4})\log_{10}{p}}$, where $r>0$ is to control signal strength and $t\in[0,1]$ is to adjust signal sparsity. In all simulations, we take $r=0.03,0.06,0.09,0.12$ and $t=0.1,0.15,0.2,0.25$. 

In this section, we only display the simulation results for sample size ${\bf n}=(20,30,45,50)$ , and the simulation results for the other two sample sizes are included in the supplementary materials

 In Simulation $1$, we compare the simulation results of our test and the three competing tests $T_{Z}$, $T_{ZZG}$, and $T_{ZZ}$ for the MANOVA problem. For the covariance matrices $\bSigma_{1}$, $\bSigma_{2}$, $\bSigma_{3}$ and $\bSigma_{4}$, we consider the following two cases respectively: 
\begin{itemize}
\item []Case 1: $\bSigma_{1}=\bSigma_{2}=\bSigma_{3}=\bSigma_{4}=3\I_{p\times p}$.

\item []Case 2: $\bSigma_{1}=\left(0.4^{|i-j|}\right)_{p\times p}$, $\bSigma_{2}=2\bSigma_{1}$, $\bSigma_{3}=1.5\bSigma_{1}$, $\bSigma_{4}=4\bSigma_{1}$.
\end{itemize} 

Tables 1 and 2 and Figures 1-6 respectively display the empirical size and power of the four tests in Simulation 1 under two different cases. From Tables 1 and 2, it can be seen that under Model 1, the size of all tests are close to the normal level $\alpha=0.05$. Under Model 2 and 3, the size of $T_{New}$ and $T_{Z}$ are close to the normal level of $0.05$, while most of the size of $T_{ZZG}$ and $T_{ZZ}$ are close to 0.04. As can be seen from Figure 1-6, the power of our test is significantly greater than that of the other three tests, and the power of $T_{Z}$, $T_{ZZG}$ and $T_{ZZ}$ is close. When we fix $t$, the power of all four tests increases as $r$ increases. When we fix $r$, as $t$ increases, the power of all tests decreases and the difference between $T_{New}$ and the other three tests becomes smaller. It shows that our test for sparse alternatives is more powerful than the other three tests.

\begin{table}[H]
\begin{center}
\textbf{Table 1}~~Empirical sizes of $T_{New}$, $T_{Z}$, $T_{ZZG}$ and $T_{ZZ}$ in Case 1.\\
\label{tab1}
\resizebox{\textwidth}{30mm}{
\begin{tabular}{cccccc|cccc|cccc}
\hline
\multirow{2}{*}{$\bf{n}$}   & \multirow{2}{*}{$p$} & \multicolumn{4}{c}{$z_{\alpha ik}\stackrel{i.i.d}\sim \mathcal{N}(0,1)$}                        & \multicolumn{4}{c}{$z_{\alpha ik}\stackrel{i.i.d}\sim {(t_4/\sqrt{2})}$} & \multicolumn{4}{c}{$z_{\alpha ik}\stackrel{i.i.d}\sim {(\chi_1^2-1)/\sqrt{2}}$} \\ \cline{3-14}
                     &                     & $T_{New}$    & $T_Z$    & $T_{ZZG}$     & \multicolumn{1}{c|}{$T_{ZZ}$}    & $T_{New}$     & $T_Z$    & $T_{ZZG}$     & \multicolumn{1}{c|}{$T_{ZZ}$}    & $T_{New}$         & $T_Z$         & $T_{ZZG}$           & $T_{ZZ}$        \\ \hline
\multirow{5}{*}{(20,30,45,50)}  & 100                  & 0.0545 &0.0575 &0.0545 &0.0510 &0.0550 &0.0575 &0.0385 &0.0325 &0.0510 &0.0555 &0.0365 &0.0335 \\     
                     & 200                  & 0.0480 &0.0470 &0.0455 &0.0445 &0.0525 &0.0540 &0.0375 &0.0335 &0.0560 &0.0605 &0.0410 &0.0365\\    
                     & 300                  & 0.0535 &0.0545 &0.0540 &0.0515 &0.0545 &0.0565 &0.0370 &0.0320 &0.0525 &0.0540 &0.0350 &0.0335 \\
                     & 400                  & 0.0530 &0.0485 &0.0490 &0.0465 &0.0550 &0.0595 &0.0380 &0.0340 &0.0495 &0.0470 &0.0340 &0.0310 \\
                     & 500                 & 0.0500 &0.0550 &0.0560 &0.0535 &0.0535 &0.0555 &0.0365 &0.0320 &0.0525 &0.0565 &0.0325 &0.0315 \\     \hline  
\end{tabular}}
\end{center}
\end{table}

\begin{table}[H]
\begin{center}
\textbf{Table 2}~~Empirical sizes of $T_{New}$, $T_{Z}$, $T_{ZZG}$ and $T_{ZZ}$ in Case 2.\\
\label{tab2}
\resizebox{\textwidth}{30mm}{
\begin{tabular}{cccccc|cccc|cccc}
\hline
\multirow{2}{*}{$\bf{n}$}   & \multirow{2}{*}{$p$} & \multicolumn{4}{c}{$z_{\alpha ik}\stackrel{i.i.d}\sim \mathcal{N}(0,1)$}                        & \multicolumn{4}{c}{$z_{\alpha ik}\stackrel{i.i.d}\sim {(t_4/\sqrt{2})}$} & \multicolumn{4}{c}{$z_{\alpha ik}\stackrel{i.i.d}\sim {(\chi_1^2-1)/\sqrt{2}}$} \\ \cline{3-14}
                     &                     & $T_{New}$    & $T_Z$    & $T_{ZZG}$     & \multicolumn{1}{c|}{$T_{ZZ}$}    & $T_{New}$     & $T_Z$    & $T_{ZZG}$     & \multicolumn{1}{c|}{$T_{ZZ}$}    & $T_{New}$         & $T_Z$         & $T_{ZZG}$           & $T_{ZZ}$        \\ \hline
\multirow{5}{*}{(20,30,45,50)}  & 100                  & 0.0505 &0.0600 &0.0590 &0.0575 &0.0515 &0.0630 &0.0430 &0.0405 &0.0575 &0.0595 &0.0440 &0.0415 \\     
                     & 200                  & 0.0575 &0.0605 &0.0590 &0.0570 &0.0510 &0.0585 &0.0480 &0.0455 &0.0555 &0.0595 &0.0455 &0.0435\\    
                     & 300                  & 0.0480 &0.0445 &0.0435 &0.0410 &0.0540 &0.0560 &0.0420 &0.0390 &0.0535 &0.0550 &0.0435 &0.0410  \\
                     & 400                  & 0.0550 &0.0570 &0.0565 &0.0540 &0.0535 &0.0510 &0.0375 &0.0350 &0.0530 &0.0565 &0.0430 &0.0400  \\
                     & 500                 & 0.0570 &0.0580 &0.0570 &0.0545 &0.0555 &0.0535 &0.0415 &0.0380 &0.0510 &0.0535 &0.0460 &0.0445  \\     \hline
\end{tabular}}
\end{center}
\end{table}

\begin{figure}[H]
\centering
    \includegraphics[scale=1]{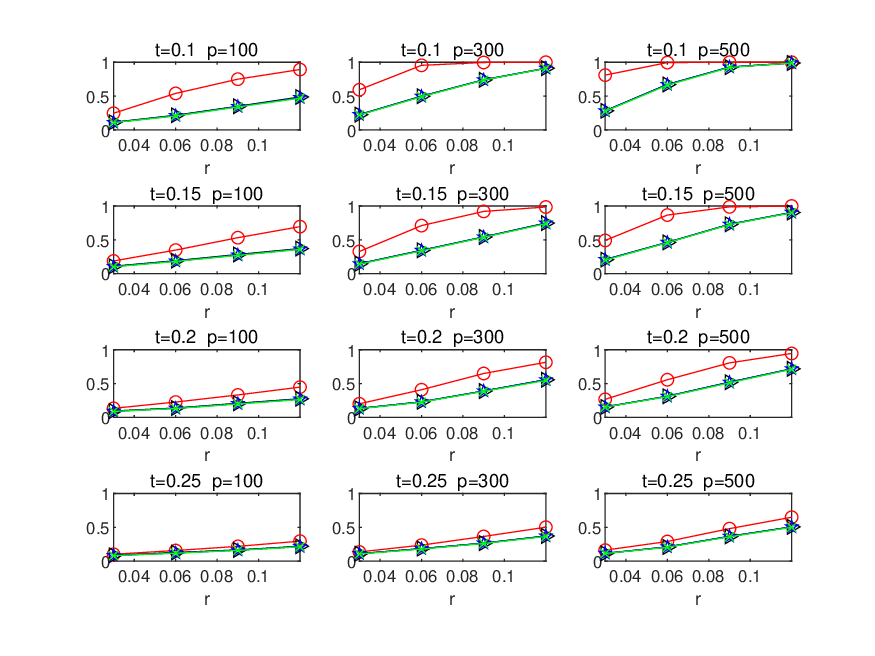}
    {\caption*{\footnotesize\textbf{Figure 1.} For Case 1, under different values of $r$ and $t$, the empirical power of the four tests when $Z_{\alpha i} \sim N(0,1)$, and $(n_{1},n_{2},n_{3},n_{4})=(20,30,45,50)$. ($T_{New}:red;  T_{Z}:black;  T_{ZZG}:blue;  T_{ZZ}:green$)}}
\end{figure}

\begin{figure}[H]
\centering
    \includegraphics[scale=1]{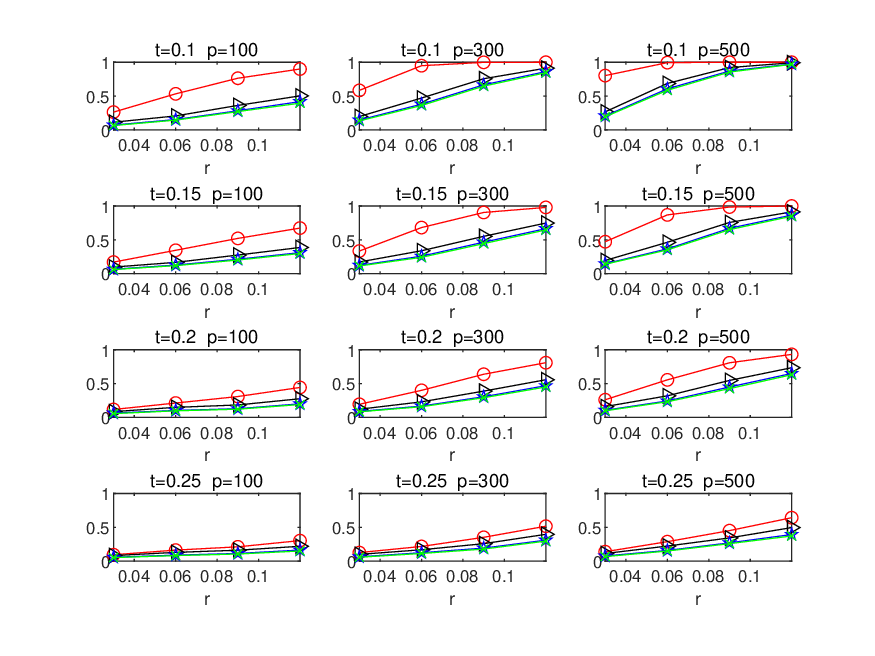}
    {\caption*{\footnotesize\textbf{Figure 2.} For Case 1, under different values of $r$ and $t$, the empirical power of the four tests when $Z_{\alpha i} \sim t_{4}/\sqrt{2}$, and $(n_{1},n_{2},n_{3},n_{4})=(20,30,45,50)$. ($T_{New}:red;  T_{Z}:black;  T_{ZZG}:blue;  T_{ZZ}:green$)}}
\end{figure}

\begin{figure}[H]
\centering
    \includegraphics[scale=1]{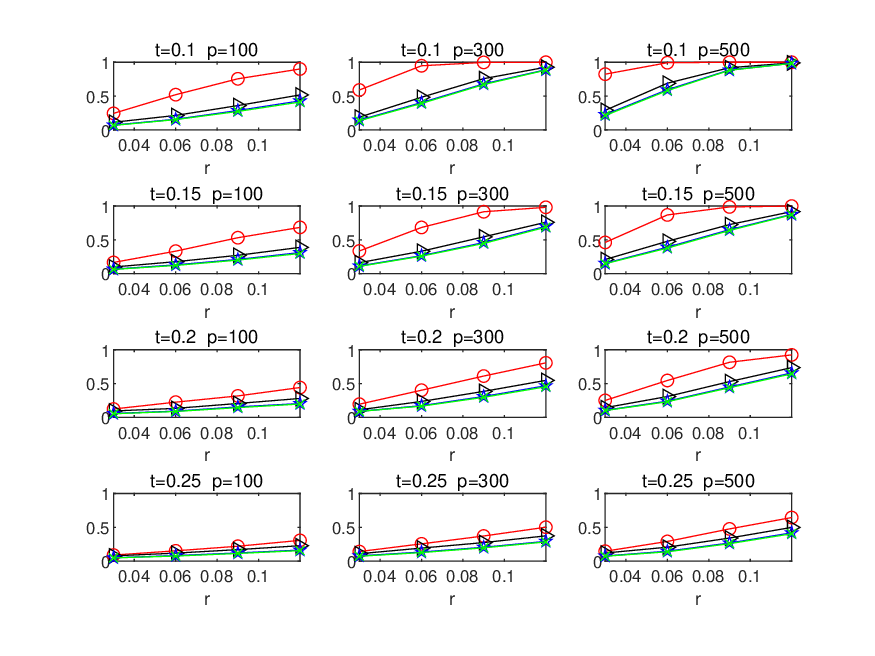}
    {\caption*{\footnotesize\textbf{Figure 3.} For Case 1, under different values of $r$ and $t$, the empirical power of the four tests when $Z_{\alpha i} \sim {(\chi_1^2-1)/\sqrt{2}}$, and $(n_{1},n_{2},n_{3},n_{4})=(20,30,45,50)$. ($T_{New}:red;  T_{Z}:black;  T_{ZZG}:blue;  T_{ZZ}:green$)}}
\end{figure}

\begin{figure}[H]
\centering
    \includegraphics[scale=1]{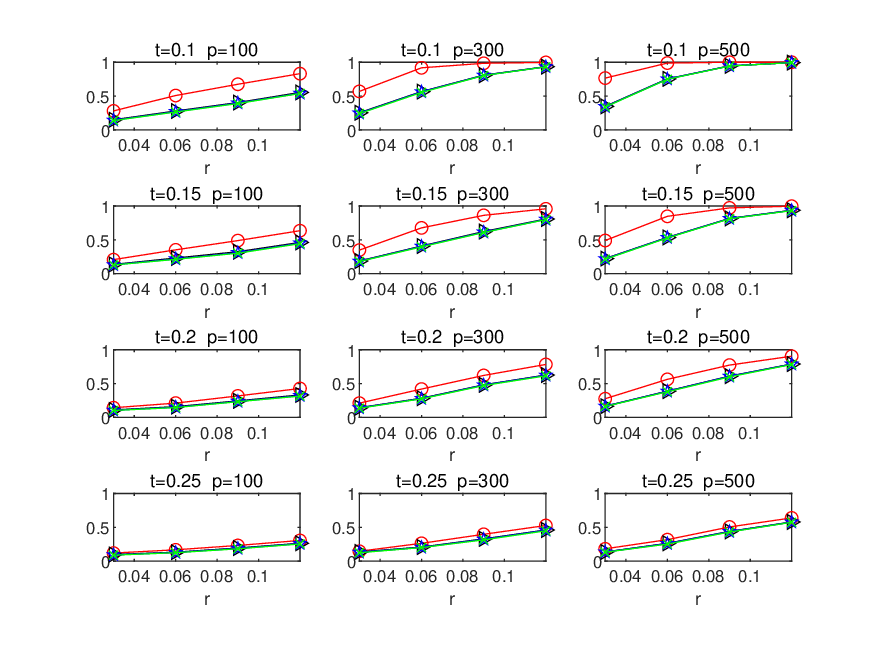}
    {\caption*{\footnotesize\textbf{Figure 4.} For Case 2, under different values of $r$ and $t$, the empirical power of the four tests when $Z_{\alpha i} \sim N(0,1)$, and $(n_{1},n_{2},n_{3},n_{4})=(20,30,45,50)$. ($T_{New}:red;  T_{Z}:black;  T_{ZZG}:blue;  T_{ZZ}:green$)}}
\end{figure}

\begin{figure}[H]
\centering
    \includegraphics[scale=1]{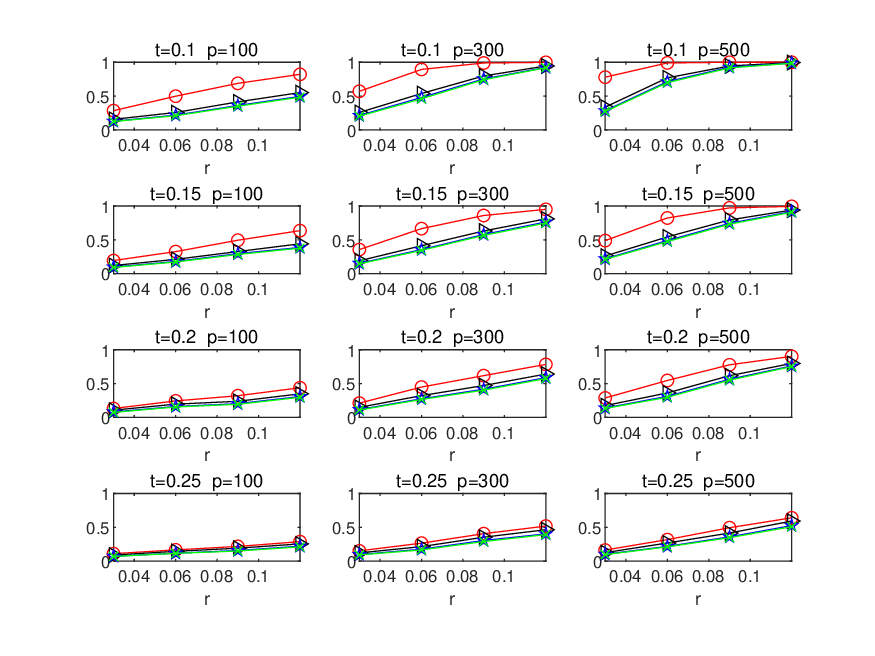}
    {\caption*{\footnotesize\textbf{Figure 5.} For Case 2, under different values of $r$ and $t$, the empirical power of the four tests when $Z_{\alpha i} \sim t_{4}/\sqrt{2}$, and $(n_{1},n_{2},n_{3},n_{4})=(20,30,45,50)$. ($T_{New}:red;  T_{Z}:black;  T_{ZZG}:blue;  T_{ZZ}:green$)}}
\end{figure}

\begin{figure}[H]
\centering
    \includegraphics[scale=1]{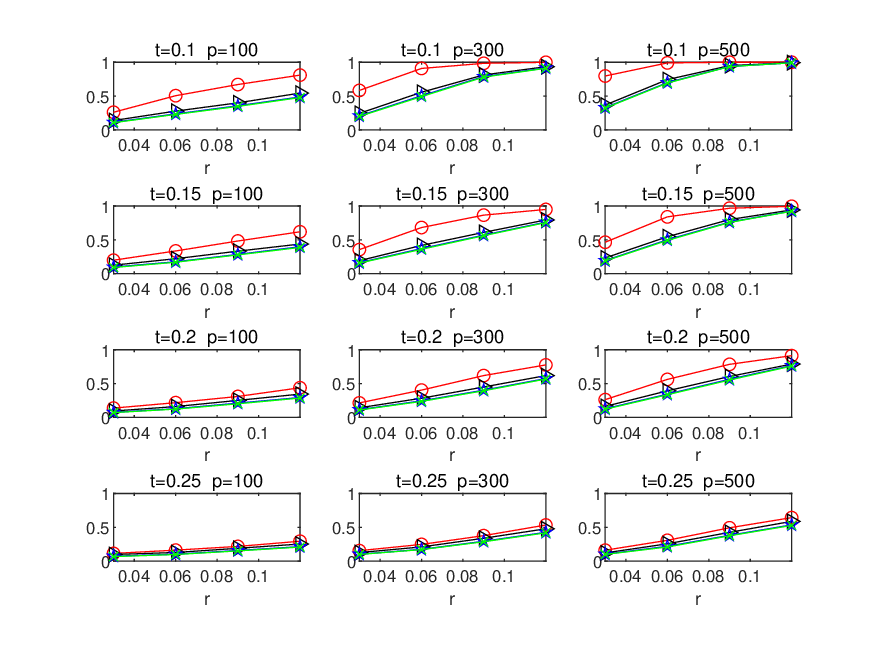}
    {\caption*{\footnotesize\textbf{Figure 6.} For Case 2, under different values of $r$ and $t$, the empirical power of the four tests when $Z_{\alpha i} \sim {(\chi_1^2-1)/\sqrt{2}}$, and $(n_{1},n_{2},n_{3},n_{4})=(20,30,45,50)$. ($T_{New}:red;  T_{Z}:black;  T_{ZZG}:blue;  T_{ZZ}:green$)}}
\end{figure}

 \clearpage

In Simulation $2$, we consider simulating the performance of four tests for GLHT problems. We chose the following special linear hypothesis: 
\begin{eqnarray}\label{eq3.1}
H_0:2\bmu_{1}-2\bmu_{2}-\bmu_{3}+3\bmu_{4}=\0.
\end{eqnarray}
Therefore, we set $\widetilde{\G}= (2\e_{1,K}-2\e_{2,K}-\e_{3,K}+3\e_{4,K})$ where $\e_{I,K}$ represents the $K$-dimensional unit vector with the $I$-th term being $1$ and the other terms being $0$. For the group covariance matrix, we have adopted the following two cases: 
\begin{itemize}
\item []Case 3: $\bSigma_{1}=\bSigma_{2}=\bSigma_{3}=\bSigma_{4}=4\I_{p\times p}$.

\item []Case 4: $\bSigma_{1}=\left(0.5^{|i-j|}\right)_{p\times p}$, $\bSigma_{2}=1.5\bSigma_{1}$, $\bSigma_{3}=2.5\bSigma_{1}$, $\bSigma_{4}=3\bSigma_{1}$.
\end{itemize}

In Case 3 and 4, the results of the four tests for Simulation 2 are shown in Table 3, 4 and Figure 7-12 respectively. According to Table 3 and 4, we can see that the size of $T_{New}$ and $T_{Z}$ is well controlled, both around 0.05. The size of $T_{ZZG}$ and $T_{ZZ}$ is better controlled in Case 4 than in Case 3. As shown in Figure 7-12, the power of all tests increases as dimension $p$ and signal strength $r$ increase, while the power of all tests decreases as $t$ increases. However, the power of $T_{New}$ is always the largest.

In all simulations, our test $T_{New}$ performed well. Furthermore, some simulation results (provided in the supplementary materials) also indicate that our testing performance is good for different sample sizes ${\bf n}$.

\begin{table}[H]
\begin{center}
\textbf{Table 3}~~Empirical sizes of $T_{New}$, $T_{Z}$, $T_{ZZG}$ and $T_{ZZ}$ in Case 3.\\
\label{tab3}
\resizebox{\textwidth}{30mm}{
\begin{tabular}{cccccc|cccc|cccc}
\hline
\multirow{2}{*}{$\bf{n}$}   & \multirow{2}{*}{$p$} & \multicolumn{4}{c}{$z_{\alpha ik}\stackrel{i.i.d}\sim \mathcal{N}(0,1)$}                        & \multicolumn{4}{c}{$z_{\alpha ik}\stackrel{i.i.d}\sim {(t_4/\sqrt{2})}$} & \multicolumn{4}{c}{$z_{\alpha ik}\stackrel{i.i.d}\sim {(\chi_1^2-1)/\sqrt{2}}$} \\ \cline{3-14}
                     &                     & $T_{New}$    & $T_Z$    & $T_{ZZG}$     & \multicolumn{1}{c|}{$T_{ZZ}$}    & $T_{New}$     & $T_Z$    & $T_{ZZG}$     & \multicolumn{1}{c|}{$T_{ZZ}$}    & $T_{New}$         & $T_Z$         & $T_{ZZG}$           & $T_{ZZ}$        \\ \hline
\multirow{5}{*}{(20,30,45,50)}  & 100                  & 0.0520 &0.0570 &0.0495 &0.0470 &0.0520 &0.0695 &0.0490 &0.0450 &0.0565 &0.0530 &0.0375 &0.0365  \\     
                     & 200                  & 0.0500 &0.0440 &0.0400 &0.0390 &0.0530 &0.0540 &0.0445 &0.0420 &0.0540 &0.0555 &0.0420 &0.0395\\    
                     & 300                  & 0.0535 &0.0495 &0.0470 &0.0445 &0.0525 &0.0515 &0.0395 &0.0360 &0.0510 &0.0555 &0.0385 &0.0370\\
                     & 400                  & 0.0525 &0.0460 &0.0420 &0.0415 &0.0505 &0.0460 &0.0330 &0.0310 &0.0520 &0.0505 &0.0380 &0.0355  \\
                     & 500                 & 0.0495 &0.0510 &0.0505 &0.0500 &0.0540 &0.0525 &0.0425 &0.0395 &0.0550 &0.0560 &0.0425 &0.0395 \\     \hline
\end{tabular}}
\end{center}
\end{table}

\begin{table}[H]
\begin{center}
\textbf{Table 4}~~Empirical sizes of $T_{New}$, $T_{Z}$, $T_{ZZG}$ and $T_{ZZ}$ in Case 4.\\
\label{tab4}
\resizebox{\textwidth}{30mm}{
\begin{tabular}{cccccc|cccc|cccc}
\hline
\multirow{2}{*}{$\bf{n}$}   & \multirow{2}{*}{$p$} & \multicolumn{4}{c}{$z_{\alpha ik}\stackrel{i.i.d}\sim \mathcal{N}(0,1)$}                        & \multicolumn{4}{c}{$z_{\alpha ik}\stackrel{i.i.d}\sim {(t_4/\sqrt{2})}$} & \multicolumn{4}{c}{$z_{\alpha ik}\stackrel{i.i.d}\sim {(\chi_1^2-1)/\sqrt{2}}$} \\ \cline{3-14}
                     &                     & $T_{New}$    & $T_Z$    & $T_{ZZG}$     & \multicolumn{1}{c|}{$T_{ZZ}$}    & $T_{New}$     & $T_Z$    & $T_{ZZG}$     & \multicolumn{1}{c|}{$T_{ZZ}$}    & $T_{New}$         & $T_Z$         & $T_{ZZG}$           & $T_{ZZ}$        \\ \hline
\multirow{5}{*}{(20,30,45,50)}  & 100                  & 0.0480 &0.0525 &0.0445 &0.0440 &0.0540 &0.0525 &0.0405 &0.0400 &0.0595 &0.0550 &0.0435 &0.0430 \\     
                     & 200                  & 0.0515 &0.0560 &0.0500 &0.0490 &0.0525 &0.0605 &0.0480 &0.0465 &0.0560 &0.0610 &0.0515 &0.0505\\    
                     & 300                  & 0.0585 &0.0590 &0.0570 &0.0565 &0.0530 &0.0545 &0.0425 &0.0415 &0.0490 &0.0515 &0.0420 &0.0400  \\
                     & 400                  & 0.0550 &0.0560 &0.0490 &0.0485 &0.0510 &0.0475 &0.0400 &0.0395 &0.0525 &0.0620 &0.0535 &0.0530  \\
                     & 500                 & 0.0555 &0.0525 &0.0490 &0.0480 &0.0505 &0.0515 &0.0430 &0.0425 &0.0505 &0.0530 &0.0460 &0.0450  \\     \hline
\end{tabular}}
\end{center}
\end{table}

\begin{figure}[H]
\centering
    \includegraphics[scale=1]{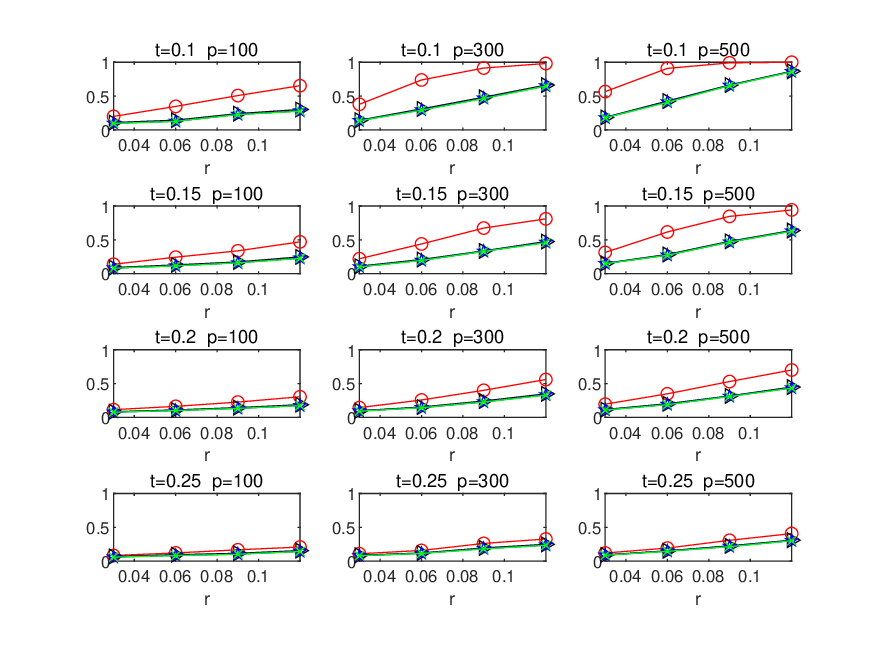}
    {\caption*{\footnotesize\textbf{Figure 7.} For Case 3, under different values of $r$ and $t$, the empirical power of the four tests when $Z_{\alpha i} \sim N(0,1)$, and $(n_{1},n_{2},n_{3},n_{4})=(20,30,45,50)$. ($T_{New}:red;  T_{Z}:black;  T_{ZZG}:blue;  T_{ZZ}:green$)}}
\end{figure}

\begin{figure}[H]
\centering
    \includegraphics[scale=1]{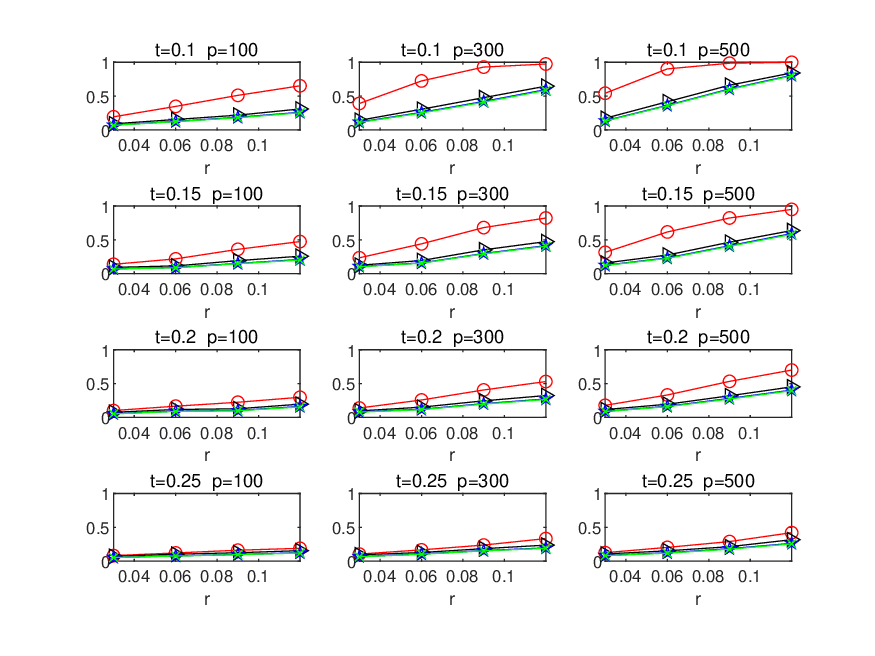}
    {\caption*{\footnotesize\textbf{Figure 8.} For Case 3, under different values of $r$ and $t$, the empirical power of the four tests when $Z_{\alpha i} \sim t_{4}/\sqrt{2}$, and $(n_{1},n_{2},n_{3},n_{4})=(20,30,45,50)$. ($T_{New}:red;  T_{Z}:black;  T_{ZZG}:blue;  T_{ZZ}:green$)}}
\end{figure}

\begin{figure}[H]
\centering
    \includegraphics[scale=1]{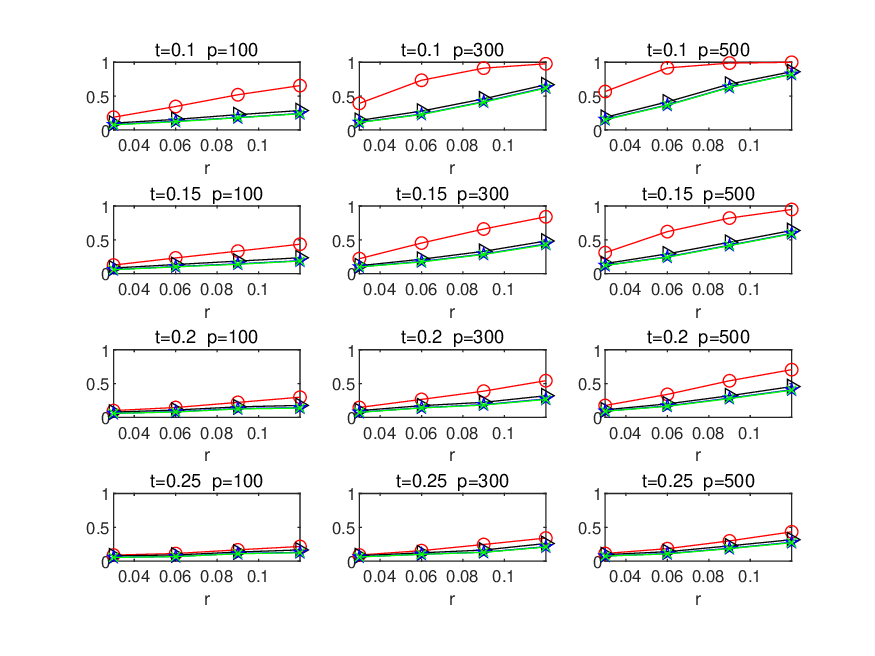}
    {\caption*{\footnotesize\textbf{Figure 9.} For Case 3, under different values of $r$ and $t$, the empirical power of the four tests when $Z_{\alpha i} \sim {(\chi_1^2-1)/\sqrt{2}}$, and $(n_{1},n_{2},n_{3},n_{4})=(20,30,45,50)$. ($T_{New}:red;  T_{Z}:black;  T_{ZZG}:blue;  T_{ZZ}:green$)}}
\end{figure}

\begin{figure}[H]
\centering
    \includegraphics[scale=1]{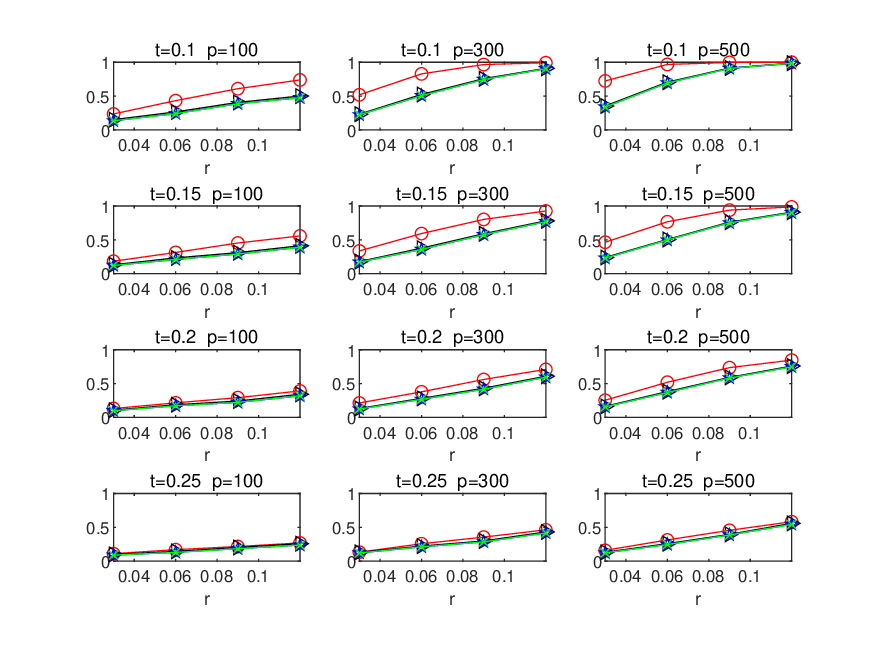}
    {\caption*{\footnotesize\textbf{Figure 10.} For Case 4, under different values of $r$ and $t$, the empirical power of the four tests when $Z_{\alpha i} \sim N(0,1)$, and $(n_{1},n_{2},n_{3},n_{4})=(20,30,45,50)$. ($T_{New}:red;  T_{Z}:black;  T_{ZZG}:blue;  T_{ZZ}:green$)}}
\end{figure}

\begin{figure}[H]
\centering
    \includegraphics[scale=1]{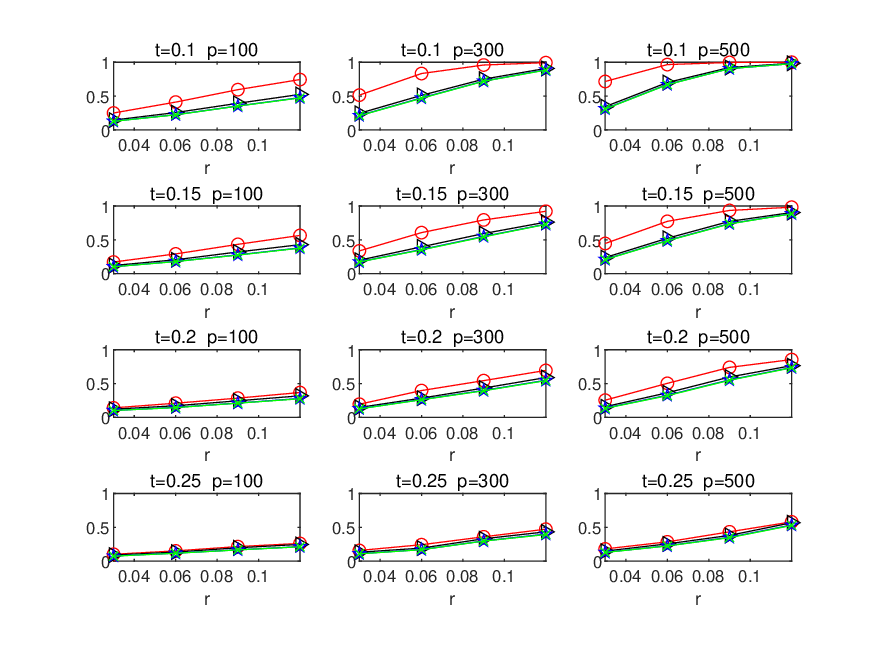}
    {\caption*{\footnotesize\textbf{Figure 11.} For Case 4, under different values of $r$ and $t$, the empirical power of the four tests when $Z_{\alpha i} \sim t_{4}/\sqrt{2}$, and $(n_{1},n_{2},n_{3},n_{4})=(20,30,45,50)$. ($T_{New}:red;  T_{Z}:black;  T_{ZZG}:blue;  T_{ZZ}:green$)}}
\end{figure}

\begin{figure}[H]
\centering
    \includegraphics[scale=1]{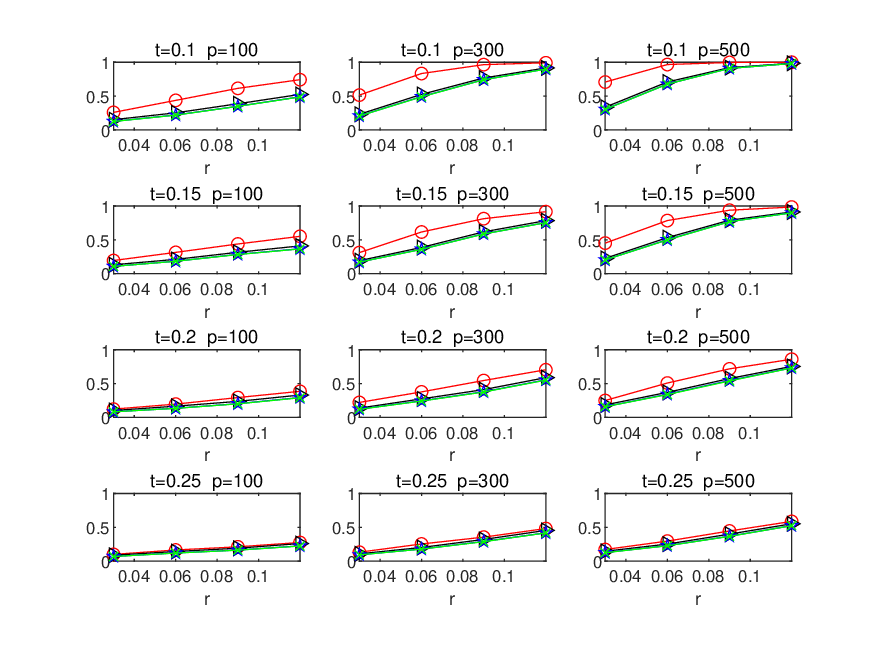}
    {\caption*{\footnotesize\textbf{Figure 12.} For Case 4, under different values of $r$ and $t$, the empirical power of the four tests when $Z_{\alpha i} \sim {(\chi_1^2-1)/\sqrt{2}}$, and $(n_{1},n_{2},n_{3},n_{4})=(20,30,45,50)$. ($T_{New}:red;  T_{Z}:black;  T_{ZZG}:blue;  T_{ZZ}:green$)}}
\end{figure}

\section{A Real Data Example}\label{sec4}

In this section, we will analyze the data of the human corneal surface. For a detailed explanation of these data, please refer to \cite{Locantore et al.:1999}. After data processing, \cite{Zhou et al.:2017} obtained 150 corneal surface datasets, each with a length of 2000. These data were divided into four groups. The normal cornea group consisted of 43 healthy corneas. The other three groups were abnormal corneas: unilateral suspect, suspect map and clinical keratoconus groups. They consisted of 14,21 and 72 corneas, respectively. Based on these datasets, \cite{Zhou et al.:2017} applied them to the GLHT problem, while \cite{Zhang and Zhu:2022b} applied it to the $k$-nearest neighbors classifier. Here, we apply these data sets to the methods of $T_{New}$, $T_{Z}$, $T_{ZZG}$ and $T_{ZZ}$ to compare the performance of our tests.

Firstly, we conduct data simulation for the one-way MANOVA problem. The simulation results are presented in Table 5. According to Table 5, the $p$-values of all four tests are far below the significance level of $0.05$, and the $p$-value of the test proposed by us is the smallest. Therefore, we have sufficient reasons to reject the null hypothesis.

\begin{table}[H]
\centering
\textbf{Table 5}~~Simulation results of corneal data set for one-way MANOVA problem\\
   \label{tab5}
   \begin{tabular}{ccccc}
      \toprule
       Method & $T_{New}$ & $T_{Z}$ & $T_{ZZG}$ & $T_{ZZ}$ \\ 
       \midrule
       Statistic & $15.8946$ & $0.8080$ & $159.7330$ & $121.1992$ \\
$p$-value & $0$ & $1.1768\times10^{-10}$ & $2.5771\times10^{-4}$ & $3.6852\times10^{-4}$ \\
d.f & --- & --- & $6.1652$ & $6.0231$\\
       \bottomrule
   \end{tabular}
\end{table}

Secondly, since the MANOVA question is rejected, it is meaningful for us to test the GLHT question. Therefore, we consider testing the following special linear hypothesis:
\begin{eqnarray}\label{eq4.1}
H_0:9\bmu_{1}-8\bmu_{2}+\bmu_{3}-\bmu_{4}=\0.
\end{eqnarray} 
The simulation results are arranged in Table 6. As shown in Table 6, all tests reject the hypothesis of \eqref{eq4.1}. In particular, the $p$-value of our test is minimal, which indicates that our proposed test procedure has sufficient evidence to reject the hypothesis.

\begin{table}[H]
\centering
\textbf{Table 6}~~Simulation results of corneal data set for linear combination of means problem.\\
   \label{tab6}
   \begin{tabular}{ccccc}
      \toprule
       Method & $T_{New}$ & $T_{Z}$ & $T_{ZZG}$ & $T_{ZZ}$ \\ 
       \midrule
       Statistic & $14.2839$ & $0.5918$ & $186.3890$ & $88.7640$ \\
$p$-value & $0$ & $5.0782\times10^{-13}$ & $1.9094\times10^{-6}$ & $0.0014$ \\
d.f & --- & --- & $2.0432$ & $1.5338$\\
       \bottomrule
   \end{tabular}
\end{table}

Since the four corneal datasets have significantly different mean corneal surfaces, we are now interested in testing whether any two corneal datasets have the same mean corneal surface. We use “Nor”, “Uni”, “Sus” and “Cli” to represent normal cornea group, unilateral suspect cornea group, suspect map cornea group, and clinical keratoconus cornea group, respectively. Table 7 shows the $p$-values of four tests under the six contrast tests. From Table 7, it can be concluded that all tests are significant for the contrast tests of Sus vs Cli and Nor vs Cli. Except for Uni vs Sus and Uni vs Nor, $T_{New}$ reject the other four contrast tests. However, $T_{Z}$ and $T_{ZZ}$ are not significant for Uni vs Sus, Uni vs Nor, Uni vs Cli and Sus vs Nor. $T_{ZZG}$ is significant for the other five contrast tests outside of Uni vs Sus, but in all six contrast tests, his approximate degrees of freedom (not shown in Table 7) are less than 3.

\begin{table}[H]
\centering
\textbf{Table 7}~~$p$-values for some contrast tests on the corneal datasets.\\
   \label{tab7}
   \begin{tabular}{ccccc}
      \toprule
      Test & $T_{New}$ & $T_{Z}$ & $T_{ZZG}$ & $T_{ZZ}$ \\ 
       \midrule
       Uni vs Sus & $0.4447$ & $0.7604$ & $0.6577$ & $0.6717$ \\
       Uni vs Nor & $0.8694$ & $0.8110$ & $0$ & $0.7527$ \\
       Uni vs Cli & $1.5508\times10^{-5}$ & $0.1591$ & $0$ & $0.1634$\\
       Sus vs Nor & $4.0601\times10^{-5}$ & $0.1052$ & $0$ & $0.1002$ \\
       Sus vs Cli & $0$ & $0$ & $0$ & $4.8453\times10^{-5}$ \\
       Nor vs Cli & $0$ & $4.4409\times10^{-16}$ & $0$ & $1.1312\times10^{-4}$ \\
       \bottomrule
   \end{tabular}
\end{table}

Finally, the example shows that our test performance is good.

\section{Concluding remarks}\label{sec5}
\noindent \quad In this article, we propose a new non parametric test based on random integration technique for the GLHT problem with high-dimensional heteroscedasticity.

Firstly, the test of high-dimensional mean vectors is an important research direction in statistics. However, due to the "curse" of the data dimension, many traditional methods have become ineffective. Therefore, there is an urgent need to develop new testing methods. We propose our own testing process by referring to the method of \cite{Jiang et al.:2024}. However, our test has a wider range of applications, not only for various post hoc and contrast tests, but also for MANOVA tests. Furthermore, our test can also be applied to high-dimensional non normal data.

Secondly, under mild conditions, we derive the asymptotic normality and power of our test. Numerical simulations show that our test is reasonable for nominal level control, and compared with some existing methods, the power of our test is superior to theirs. Through a real examples, it also shows that the performance of our test is good.

Finally, the test we propose uses the $p$-dimensional independent density function as the weight function. Therefore, we can further consider other function options. This is an interesting question that we will discuss further in the future.

\section*{Acknowledgments}

Dr Cao's research is supported by Humanities and Social Sciences Fund of the Ministry of Education (No. 22YJC910001).

\begin{appendices}
\section{Technical details}\label{appA.}
\subsection*{Proof of Theorem 2.1}
\noindent let's $\Y_1=(\Y_{11}^T,\ldots,Y_{1K}^T)^T$ be independent copies of $\Y=(Y_1^T,\ldots,\Y_K^T)^T$.
\begin{align*}
{\rm{RI}}_{w}&=\int{\E^{2}[(\G\otimes\delta^T)\Y]}\omega(\delta)d\delta\\
&=\int\E[(\G\otimes\delta^T)\Y]E[(\G\otimes\delta^T)\Y_1]\omega(\delta)d\delta\\
&=\int\E\{[(\G\otimes\delta^T)\Y][(\G\otimes\delta^T)\Y_1]\}\omega(\delta)d\delta.                                                    
\end{align*}
Then, by Fubini's theorem, we have
\begin{align*}
{\rm{RI}}&:=\int\E\{[(\G\otimes\delta^T)\Y][(\G\otimes\delta^T)\Y_1]\}\omega(\delta)d\delta  \\
&=\E\left\{\int\left[\Y_1^T(\G^T\G\otimes\delta\delta^T)\Y\right]\omega(\delta)d\delta\right\}\\
&=\E\left\{\int\left[\sum_{\alpha=1}^{K}\Y_{1\alpha}^T(\bf g_\alpha^T\bf g_\alpha\otimes\delta\delta^T)\Y_\alpha
+\sum_{\alpha\neq\beta}\Y_{1\alpha}^T(\bf g_\alpha^T\bf g_\beta\otimes\delta\delta^T)\Y_\beta\right]\omega(\delta)d\delta\right\}\\ 
&=\E\left\{\left[\sum_{\alpha=1}^{K}d_{\alpha\alpha}\Y_{1\alpha}^T(\B+\bf a\bf a^T)\Y_\alpha
+\sum_{\alpha\neq\beta}d_{\alpha\beta}\Y_{1\alpha}^T(\B+\bf a\bf a^T)\Y_\beta\right]\right\}\\
&=\left\{\left[\sum_{\alpha=1}^{K}d_{\alpha\alpha}\bmu_{\alpha}^T(\B+\bf a\bf a^T)\bmu_\alpha
+\sum_{\alpha\neq\beta}d_{\alpha\beta}\bmu_{\alpha}^T(\B+\bf a\bf a^T)\bmu_\beta\right]\right\}\\
&=\sum_{\alpha=1}^{K}d_{\alpha\alpha}\bmu_{\alpha}^T\W\bmu_\alpha
+\sum_{\alpha\neq\beta}d_{\alpha\beta}\bmu_{\alpha}^T\W\bmu_\beta\\
&={\lVert(\G\otimes\W^{\frac{1}{2}})\bmu \rVert}^2 ,                                                                                                     
\end{align*}
where $\W=\B+{\bf a}{\bf a}^T, {\bf a}=(\alpha_1,\ldots,\alpha_p)^T, \B=diag(\beta_1^2,\ldots,\beta_p^2)$, $d_{\alpha\beta}={\bf g}_\alpha^T{\bf g}_\beta$. This completes the proof.

\subsection*{Proof of Lemma 2.1}
\noindent It is clear from \eqref{eq2.12} and \eqref{eq2.16} that $\S_{nm}$ is square integrable martingale with zero mean. And it is known from  Assumption (A2) that $\S_{nm}$ has finite expectations. Therefore, we only need to prove that $\E(\S_{n,m+1}|\mathcal{F}_{nm})=\S_{nm}$, and this condition expectation can be written as 
$$\E(\S_{n,m+1}|\mathcal{F}_{nm})=\E(\sum_{j=2}^{m+1}\V_{nj}|\mathcal{F}_{nm})=\E(\sum_{j=2}^{m}\V_{nj}|\mathcal{F}_{mn})+\E(\V_{n,m+1}|\mathcal{F}_{nm}).$$
Since $\V_{n2},\ldots,\V_{nm}\in \mathcal{F}_{nm}$,
$$\E(\sum_{j=2}^{m}\V_{nj}|\mathcal{F}_{nm})=\S_{nm}.$$
Furthermore, since $\{\Y_{1},\ldots,\Y_{m}\}\in \mathcal{F}_{nm}$ and $\Y_{m+1}$ is independent of $\mathcal{F}_{nm}$, 
$$\E(\Y_{m+1}|\mathcal{F}_{nm})=0.$$
Then by \eqref{eq2.19}, we get
$$\E(\V_{n,m+1}|\mathcal{F}_{nm})=\E(\sum_{i=1}^{m}\phi_{i,m+1}|\mathcal{F}_{nm})=0.$$
Therefore, this completes the proof of this lemma.

\subsection*{Proof of Lemma 2.2}
\noindent For $i\in\{m_{\gamma-1}+1,\ldots,m_{\gamma}\}$, let 
$$\Y_{i}=\W^{1/2}\X_{\gamma(i-m_{\gamma-1}})=\W^{1/2}\X_{\gamma i_{0}}=\W^{1/2}\bGam_{\gamma}\Z_{\gamma i_{0}},$$ 
where $i_{0}=i-m_{\gamma-1}$, then we have 
\begin{align*}
&\E(\Y_{i}^{T}\W^{1/2}\bSigma_{\alpha}\W^{1/2}\Y_{i}\Y_{i}^{T}\W^{1/2}\bSigma_{\beta}\W^{1/2}\Y_{i})\\
&=\E(\Z_{\gamma i_{0}}^{T}\bGam_{\gamma}^{T}\W\bSigma_{\alpha}\W\bGam_{\gamma}\Z_{\gamma i_{0}}\Z_{\gamma i_{0}}^{T}\bGam_{\gamma}^{T}\W\bSigma_{\beta}\W\bGam_{\gamma}\Z_{\gamma i_{0}}).
\end{align*}
Using Proposition A.1.(i) of \cite{Chen et al.:2010}, we have
\begin{align*}
\E&(\Z_{\gamma i_{0}}^{T}\bGam_{\gamma}^{T}\W\bSigma_{\alpha}\W\bGam_{\gamma}\Z_{\gamma i_{0}}\Z_{\gamma i_{0}}^{T}\bGam_{\gamma}^{T}\W\bSigma_{\beta}\W\bGam_{\gamma}\Z_{\gamma i_{0}})\\
=&\rm{tr}(\bGam_{\gamma}^{T}\W\bSigma_{\alpha}\W\bGam_{\gamma})\rm{tr}(\bGam_{\gamma}^{T}\W\bSigma_{\beta}\W\bGam_{\gamma}+
2\rm{tr}(\bGam_{\gamma}^{T}\W\bSigma_{\alpha}\W\bGam_{\gamma}\bGam_{\gamma}^{T}\W\bSigma_{\beta}\W\bGam_{\gamma})\\
&+\Delta\rm{tr}(\bGam_{\gamma}^{T}\W\bSigma_{\alpha}\W\bGam_{\gamma}\circ\bGam_{\gamma}^{T}\W\bSigma_{\beta}\W\bGam_{\gamma})\\
\leq &\rm{tr}(\W\bSigma_{\alpha}\W\bSigma_{\gamma})\rm{tr}(\W\bSigma_{\beta}\W\bSigma_{\gamma})+
(2+\Delta)\rm{tr}(\W\bSigma_{\alpha}\W\bSigma_{\gamma}\W\bSigma_{\beta}\W\bSigma_{\gamma}).
\end{align*}

To prove the second statement, for $i\in\{m_{\alpha-1}+1,\ldots,m_{\alpha}\}$, $j\in\{m_{\beta-1}+1,\ldots,m_{\beta}\}$,  $k\in\{m_{\gamma-1}+1,\ldots,m_{\gamma}\}$, where $i\neq j$, $k\neq i$, and $k\neq j$, we denote $\Y_{i}=\W^{1/2}\bGam_{\alpha}\Z_{\alpha i_{0}}$, $\Y_{j}=\W^{1/2}\bGam_{\beta}\Z_{\beta j_{0}}$ and $\Y_{k}=\W^{1/2}\bGam_{\gamma}\Z_{\gamma k_{0}}$, where $i_{0}=i-m_{\alpha-1}$, $j_{0}=j-m_{\beta-1}$, and $k_{0}=k-m_{\gamma-1}$. So we have 
\begin{align*}
&\E\left((\Y_{i}^{T}\Y_{k})^{2}(\Y_{j}^{T}\Y_{k})^{2}\right)\\
&=\E\left(\Y_{k}^{T}\Y_{i}\Y_{i}^{T}\Y_{k}\Y_{k}^{T}\Y_{j}\Y_{j}^{T}\Y_{k}\right)\\
&=\E\left\{\E\left(\Y_{k}^{T}\Y_{i}\Y_{i}^{T}\Y_{k}\Y_{k}^{T}\Y_{j}\Y_{j}^{T}\Y_{k}|\Y_{i},\Y_{j}\right)\right\}\\
&=\E\left\{\E\left(\Z_{\gamma k_{0}}^{T}\bGam_{\gamma}^{T}\W^{1/2}\Y_{i}\Y_{i}^{T}\W^{1/2}\bGam_{\gamma}\Z_{\gamma k_{0}}\Z_{\gamma k_{0}}^{T}\bGam_{\gamma}^{T}\W^{1/2}\Y_{j}\Y_{j}^{T}\W^{1/2}\bGam_{\gamma}\Z_{\gamma k_{0}}|\Y_{i},\Y_{j}\right)\right\}.
\end{align*}
Using Proposition A.1.(i) of \cite{Chen et al.:2010}, we have
\begin{align*}
&\E\left\{\E\left(\Z_{\gamma k_{0}}^{T}\bGam_{\gamma}^{T}\W^{1/2}\Y_{i}\Y_{i}^{T}\W^{1/2}\bGam_{\gamma}\Z_{\gamma k_{0}}\Z_{\gamma 
k_{0}}^{T}\bGam_{\gamma}^{T}\W^{1/2}\Y_{j}\Y_{j}^{T}\W^{1/2}\bGam_{\gamma}\Z_{\gamma k_{0}}|\Y_{i},\Y_{j}\right)\right\}\\
&=\E\bigg\{\rm{tr}(\bGam_{\gamma}^{T}\W^{1/2}\Y_{i}\Y_{i}^{T}\W^{1/2}\bGam_{\gamma})\rm{tr}(\bGam_{\gamma}^{T}\W^{1/2}\Y_{j}\Y_{j}^{T}\W^{1/2}\bGam_{\gamma})\\
&\quad+2\rm{tr}(\bGam_{\gamma}^{T}\W^{1/2}\Y_{i}\Y_{i}^{T}\W^{1/2}\bGam_{\gamma}\bGam_{\gamma}^{T}\W^{1/2}\Y_{j}\Y_{j}^{T}\W^{1/2}\bGam_{\gamma})\\
&\quad+\Delta\rm{tr}(\bGam_{\gamma}^{T}\W^{1/2}\Y_{i}\Y_{i}^{T}\W^{1/2}\bGam_{\gamma}\circ\bGam_{\gamma}^{T}\W^{1/2}\Y_{j}\Y_{j}^{T}\W^{1/2}\bGam_{\gamma})\bigg\}\\
&\leq \E\bigg\{\rm{tr}(\bGam_{\gamma}^{T}\W^{1/2}\Y_{i}\Y_{i}^{T}\W^{1/2}\bGam_{\gamma})\rm{tr}(\bGam_{\gamma}^{T}\W^{1/2}\Y_{j}\Y_{j}^{T}\W^{1/2}\bGam_{\gamma})\\
&\quad+(2+\Delta)\rm{tr}(\bGam_{\gamma}^{T}\W^{1/2}\Y_{i}\Y_{i}^{T}\W^{1/2}\bGam_{\gamma}\bGam_{\gamma}^{T}\W^{1/2}\Y_{j}\Y_{j}^{T}\W^{1/2}\bGam_{\gamma})\bigg\}\\
&=\E\bigg\{\rm{tr}(\Y_{i}^{T}\W^{1/2}\bGam_{\gamma}\bGam_{\gamma}^{T}\W^{1/2}\Y_{i})\rm{tr}(\Y_{j}^{T}\W^{1/2}\bGam_{\gamma}\bGam_{\gamma}^{T}\W^{1/2}\Y_{j})\\
&\quad+(2+\Delta)\rm{tr}(\Y_{i}^{T}\W^{1/2}\bGam_{\gamma}\bGam_{\gamma}^{T}\W^{1/2}\Y_{j}\Y_{j}^{T}\W^{1/2}\bGam_{\gamma}\bGam_{\gamma}^{T}\W^{1/2}\Y_{i})\bigg\}\\
&=\E\bigg\{\rm{tr}(\Z_{\alpha i_{0}}^{T}\bGam_{\alpha}^{T}\W\bSigma_{\gamma}\W\bGam_{\alpha}\Z_{\alpha i_{0}})\rm{tr}(\Z_{\beta j_{0}}^{T}\bGam_{\beta}^{T}\W\bSigma_{\gamma}\W\bGam_{\beta}\Z_{\beta j_{0}})\\
&\quad+(2+\Delta)\rm{tr}(\Z_{\alpha i_{0}}^{T}\bGam_{\alpha}^{T}\W\bSigma_{\gamma}\W\bGam_{\beta}\Z_{\beta j_{0}}\Z_{\beta j_{0}}^{T}\bGam_{\beta}^{T}\W\bSigma_{\gamma}\W\bGam_{\alpha}\Z_{\alpha i_{0}})\bigg\}\\
&=\rm{tr}(\W\bSigma_{\alpha}\W\bSigma_{\gamma})\rm{tr}(\W\bSigma_{\beta}\W\bSigma_{\gamma})+
(2+\Delta)\rm{tr}(\W\bSigma_{\alpha}\W\bSigma_{\gamma}\W\bSigma_{\beta}\W\bSigma_{\gamma}).
\end{align*}
So in summary, the lemma proof is complete.

\subsection*{Proof of Lemma 2.3}
\noindent For $i\in\{m_{\alpha-1}+1,\ldots,m_{\alpha}\}$, $j\in\{m_{\beta-1}+1,\ldots,m_{\beta}\}$, we denote $\Y_{i}=\W^{1/2}\bGam_{\alpha}\Z_{\alpha i_{0}}$, $\Y_{j}=\W^{1/2}\bGam_{\beta}\Z_{\beta j_{0}}$, where $i_{0}=i-m_{\alpha-1}$, $j_{0}=j-m_{\beta-1}$. So we have
\begin{align*}
&\E(\Y_{i}^{T}\Y_{j})^{4}\\
&=\E\left\{\E\left(\Y_{i}^{T}\Y_{j}\Y_{j}^{T}\Y_{i}\Y_{i}^{T}\Y_{j}\Y_{j}^{T}\Y_{i}|\Y_{j}\right)\right\}\\
&=\E\left\{\E\left(\Z_{\alpha i_{0}}^{T}\bGam_{\alpha}^{T}\W^{1/2}\Y_{j}\Y_{j}^{T}\W^{1/2}\bGam_{\alpha}\Z_{\alpha i_{0}}\Z_{\alpha 
i_{0}}^{T}\bGam_{\alpha}^{T}\W^{1/2}\Y_{j}\Y_{j}^{T}\W^{1/2}\bGam_{\alpha}\Z_{\alpha i_{0}}|\Y_{j}\right)\right\}.
\end{align*}
Using Proposition A.1.(i) of \cite{Chen et al.:2010}, we have
\begin{align*}
&\E\left\{\E\left(\Z_{\alpha i_{0}}^{T}\bGam_{\alpha}^{T}\W^{1/2}\Y_{j}\Y_{j}^{T}\W^{1/2}\bGam_{\alpha}\Z_{\alpha i_{0}}\Z_{\alpha 
i_{0}}^{T}\bGam_{\alpha}^{T}\W^{1/2}\Y_{j}\Y_{j}^{T}\W^{1/2}\bGam_{\alpha}\Z_{\alpha i_{0}}|\Y_{j}\right)\right\}\\
&=\E\bigg\{\rm{tr}^{2}(\bGam_{\alpha}^{T}\W^{1/2}\Y_{j}\Y_{j}^{T}\W^{1/2}\bGam_{\alpha})
+2\rm{tr}(\bGam_{\alpha}^{T}\W^{1/2}\Y_{j}\Y_{j}^{T}\W^{1/2}\bGam_{\alpha}\bGam_{\alpha}^{T}\W^{1/2}\Y_{j}\Y_{j}^{T}\W^{1/2}\bGam_{\alpha})\\
&\quad+\Delta\rm{tr}(\bGam_{\alpha}^{T}\W^{1/2}\Y_{j}\Y_{j}^{T}\W^{1/2}\bGam_{\alpha}\circ\bGam_{\alpha}^{T}\W^{1/2}\Y_{j}\Y_{j}^{T}\W^{1/2}\bGam_{\alpha})\bigg\}\\
&\leq \E\bigg\{\rm{tr}^{2}(\Y_{j}^{T}\W^{1/2}\bGam_{\alpha}\bGam_{\alpha}^{T}\W^{1/2}\Y_{j})\\
&\quad+(2+\Delta)\rm{tr}(\Y_{j}^{T}\W^{1/2}\bGam_{\alpha}\bGam_{\alpha}^{T}\W^{1/2}\Y_{j}\Y_{j}^{T}\W^{1/2}\bGam_{\alpha}\bGam_{\alpha}^{T}\W^{1/2}\Y_{j})\bigg\}\\
&=\E\bigg\{(3+\Delta)(\Y_{j}^{T}\W^{1/2}\bSigma_{\alpha}\W^{1/2}\Y_{j})^{2}\bigg\}.
\end{align*}
Using Proposition A.1.(i) of \cite{Chen et al.:2010} again, we get
\begin{align*}
&\E\bigg\{(3+\Delta)(\Y_{j}^{T}\W^{1/2}\bSigma_{\alpha}\W^{1/2}\Y_{j})^{2}\bigg\}\\
&=(3+\Delta)\E\bigg\{(\Z_{\beta j_{0}}^{T}\bGam_{\beta}^{T}\W\bSigma_{\alpha}\W\bGam_{\beta}\Z_{\beta j_{0}}\Z_{\beta j_{0}}^{T}\bGam_{\beta}^{T}\W\bSigma_{\alpha}\W\bGam_{\beta}\Z_{\beta j_{0}})\bigg\}\\
&=(3+\Delta)\bigg\{\rm{tr}^{2}(\bGam_{\beta}^{T}\W\bSigma_{\alpha}\W\bGam_{\beta})\\
&\quad+2\rm{tr}(\bGam_{\beta}^{T}\W\bSigma_{\alpha}\W\bGam_{\beta}\bGam_{\beta}^{T}\W\bSigma_{\alpha}\W\bGam_{\beta})\\
&\quad+\Delta\rm{tr}(\bGam_{\beta}^{T}\W\bSigma_{\alpha}\W\bGam_{\beta}\circ\bGam_{\beta}^{T}\W\bSigma_{\alpha}\W\bGam_{\beta}\bigg\}\\
&\leq (3+\Delta)\bigg\{\rm{tr}^{2}(\bGam_{\beta}^{T}\W\bSigma_{\alpha}\W\bGam_{\beta})\\
&\quad+(2+\Delta)\rm{tr}(\bGam_{\beta}^{T}\W\bSigma_{\alpha}\W\bGam_{\beta}\bGam_{\beta}^{T}\W\bSigma_{\alpha}\W\bGam_{\beta})\bigg\}\\
&=(3+\Delta)\rm{tr}^{2}(\W\bSigma_{\alpha}\W\bSigma_{\beta})+(3+\Delta)(2+\Delta)\rm{tr}(\W\bSigma_{\alpha}\W\bSigma_{\beta})^{2}.
\end{align*}
Then the lemma proof is complete.

\subsection*{Proof of Lemma 2.4}
\noindent First, set $j\in\{m_{\beta-1}+1,\ldots,m_{\beta}\}$, we have 
\begin{align}\label{A.1}
\V_{nj}&=\sum_{i=1}^{j-1}\phi_{ij}\nonumber\\
&=\sum_{i=1}^{m_{\beta-1}}\phi_{ij}+\sum_{i=m_{\beta-1}+1}^{j-1}\phi_{ij}\nonumber\\
&=2\sum_{\alpha=1}^{\beta-1}\sum_{i=m_{\alpha-1}+1}^{m_{\alpha}}\frac{d_{\alpha\beta}}{n_{\alpha}n_{\beta}}\Y_{i}^{T}\Y_{j}
+2\sum_{i=m_{\beta-1}+1}^{j-1}\frac{d_{\beta\beta}}{n_{\beta}(n_{\beta}-1)}\Y_{i}^{T}\Y_{j}.
\end{align}
Then we have $\E(\V_{nj}|\mathcal{F}_{n,j-1})=0$ and
\begin{align*}
\V_{nj}^{2}&=4\left(\sum_{\alpha=1}^{\beta-1}\sum_{i=m_{\alpha-1}+1}^{m_{\alpha}}\sum_{\alpha_{1}=1}^{\beta-1}\sum_{i_{1}=m_{\alpha_{1}-1}+1}^{m_{\alpha_{1}}}
\frac{d_{\alpha\beta}d_{\alpha_{1}\beta}}{n_{\alpha}n_{\alpha_{1}}n_{\beta}^{2}}\Y_{j}^{T}\Y_{i}\Y_{i_{1}}^{T}\Y_{j}\right)\\
&\quad+4\left(\sum_{i=m_{\beta-1}+1}^{j-1}\sum_{i_{1}=m_{\beta-1}+1}^{j-1}\frac{d_{\beta\beta}^{2}}{n_{\beta}^{2}(n_{\beta}-1)^{2}}\Y_{j}^{T}\Y_{i}\Y_{i_{1}}^{T}\Y_{j}\right)\\
&\quad+8\left(\sum_{\alpha=1}^{\beta-1}\sum_{i=m_{\alpha-1}+1}^{m_{\alpha}}\sum_{i_{1}=m_{\beta-1}+1}^{j-1}
\frac{d_{\alpha\beta}d_{\beta\beta}}{n_{\alpha}n_{\beta}^{2}(n_{\beta}-1)}\Y_{j}^{T}\Y_{i}\Y_{i_{1}}^{T}\Y_{j}\right).
\end{align*}
So, 
\begin{align}\label{A.2}
\E(\V_{nj}^{2}|\mathcal{F}_{n,j-1})&=4\left(\sum_{\alpha=1}^{\beta-1}\sum_{i=m_{\alpha-1}+1}^{m_{\alpha}}\sum_{\alpha_{1}=1}^{\beta-1}\sum_{i_{1}=m_{\alpha_{1}-1}+1}^{m_{\alpha_{1}}}
\frac{d_{\alpha\beta}d_{\alpha_{1}\beta}}{n_{\alpha}n_{\alpha_{1}}n_{\beta}^{2}}\Y_{i}^{T}\W^{1/2}\bSigma_{\beta}\W^{1/2}\Y_{i_{1}}\right)\nonumber\\
&\quad+4\left(\sum_{i=m_{\beta-1}+1}^{j-1}\sum_{i_{1}=m_{\beta-1}+1}^{j-1}\frac{d_{\beta\beta}^{2}}{n_{\beta}^{2}(n_{\beta}-1)^{2}}\Y_{i}^{T}\W^{1/2}\bSigma_{\beta}\W^{1/2}\Y_{i_{1}}\right)\nonumber\\
&\quad+8\left(\sum_{\alpha=1}^{\beta-1}\sum_{i=m_{\alpha-1}+1}^{m_{\alpha}}\sum_{i_{1}=m_{\beta-1}+1}^{j-1}
\frac{d_{\alpha\beta}d_{\beta\beta}}{n_{\alpha}n_{\beta}^{2}(n_{\beta}-1)}\Y_{i}^{T}\W^{1/2}\bSigma_{\beta}\W^{1/2}\Y_{i_{1}}\right)\nonumber\\
&:=4I_{1}+4I_{2}+8I_{3},
\end{align}
and 
\begin{align}\label{A.3}
\E(\V_{nj}^{2})&=4\left(\sum_{\alpha=1}^{\beta-1}\frac{d_{\alpha\beta}^{2}}{n_{\alpha}n_{\beta}^{2}}tr(\W\bSigma_{\alpha}\W\bSigma_{\beta})\right)
+4\left(\frac{(j-1-m_{\beta-1})d_{\beta\beta}^{2}}{n_{\beta}^{2}(n_{\beta}-1)^{2}}tr(\W\bSigma_{\beta})^{2}\right).
\end{align}
Thus,
\begin{align*}
\E\left(\sum_{j=1}^{n}\E(\V_{nj}^{2}|\mathcal{F}_{n,j-1})\right)=\sum_{j=1}^{n}\E(\V_{nj}^{2})
=\sum_{\beta=1}^{K}\sum_{j=m_{\beta-1}+1}^{m_{\beta}}\E(\V_{nj}^{2})=\sigma^{2}({\rm{T}}_{n0}).
\end{align*}
Further, from \eqref{A.2}, we show by calculation that 
\begin{align*}
&\E\left[\sum_{j=m_{\beta-1}+1}^{m_{\beta}}\E(\V_{nj}^{2}|\mathcal{F}_{n,j-1})\right]^{2}\\
&=\E\left[\sum_{j=m_{\beta-1}+1}^{m_{\beta}}(4I_{1}+4I_{2}+8I_{3})\sum_{j_{1}=m_{\beta-1}+1}^{m_{\beta}}(4I_{1}^{'}+4I_{2}^{'}+8I_{3}^{'})\right]\\
&=\E\left[\sum_{j=m_{\beta-1}+1}^{m_{\beta}}(16I_{1}^{2}+16I_{2}^{2}+64I_{3}^{2}+32I_{1}I_{2})\right]\\
&\quad+\E\left[2\sum_{j=m_{\beta-1}+1}^{m_{\beta}}\sum_{j<j_{1}}^{m_{\beta}}(16I_{1}I_{1}^{'}+16I_{2}I_{2}^{'}+64I_{3}I_{3}^{'}+16I_{1}I_{2}^{'}+16I_{2}I_{1}^{'})\right]\\
&=\sum_{j=m_{\beta-1}+1}^{m_{\beta}}[16\E(I_{1}^{2})+16\E(I_{2}^{2})+64\E(I_{3}^{2})+32\E(I_{1}I_{2})]\\
&\quad+2\sum_{j=m_{\beta-1}+1}^{m_{\beta}}\sum_{j<j_{1}}^{m_{\beta}}[16\E(I_{1}I_{1}^{'})+16\E(I_{2}I_{2}^{'})+64\E(I_{3}I_{3}^{'})+16\E(I_{1}I_{2}^{'})+16\E(I_{2}I_{1}^{'})].
\end{align*}

The term $\E(I_{1}^{2})$ can be split into $\E(I_{1}^{2})=I_{11}+I_{121}+I_{122}$. This representation is the result of a combination of the following five cases:\\
(1)$\quad(\alpha=\alpha_{1}=\alpha_{2}=\alpha_{3}; i=i_{1}=i_{2}=i_{3})$
\begin{align*}
\E(I_{1}^{2})&=\E\left\{\sum_{\alpha=1}^{\beta-1}\sum_{i=m_{\alpha-1}+1}^{m_{\alpha}}
\frac{d_{\alpha\beta}^{4}}{n_{\alpha}^{4}n_{\beta}^{4}}(\Y_{i}^{T}\W^{1/2}\bSigma_{\beta}\W^{1/2}\Y_{i})^{2}\right\}\\
&=\sum_{\alpha=1}^{\beta-1}\sum_{i=m_{\alpha-1}+1}^{m_{\alpha}}
\frac{d_{\alpha\beta}^{4}}{n_{\alpha}^{4}n_{\beta}^{4}}\E(\Y_{i}^{T}\W^{1/2}\bSigma_{\beta}\W^{1/2}\Y_{i})^{2},
\end{align*}
(2)$\quad(\alpha=\alpha_{1}=\alpha_{2}=\alpha_{3}; i=i_{1}\neq i_{2}=i_{3})$
\begin{align*}
\E(I_{1}^{2})&=\E\left\{\sum_{\alpha=1}^{\beta-1}\sum_{i\neq i_{2}}^{m_{\alpha}}
\frac{d_{\alpha\beta}^{4}}{n_{\alpha}^{4}n_{\beta}^{4}}(\Y_{i}^{T}\W^{1/2}\bSigma_{\beta}\W^{1/2}\Y_{i}\Y_{i_{2}}^{T}\W^{1/2}\bSigma_{\beta}\W^{1/2}\Y_{i_{2}})\right\}\\
&=\sum_{\alpha=1}^{\beta-1}\sum_{i\neq i_{2}}^{m_{\alpha}}
\frac{d_{\alpha\beta}^{4}}{n_{\alpha}^{4}n_{\beta}^{4}}\rm{tr}^{2}(\W\bSigma_{\alpha}\W\bSigma_{\beta}),
\end{align*}
(3)$\quad(\alpha=\alpha_{1}=\alpha_{2}=\alpha_{3}; i=i_{2}\neq i_{1}=i_{3})$ or $(\alpha=\alpha_{1}=\alpha_{2}=\alpha_{3}; i=i_{3}\neq i_{1}=i_{2})$
\begin{align*}
\E(I_{1}^{2})&=2\E\left\{\sum_{\alpha=1}^{\beta-1}\sum_{i\neq i_{1}}^{m_{\alpha}}
\frac{d_{\alpha\beta}^{4}}{n_{\alpha}^{4}n_{\beta}^{4}}(\Y_{i}^{T}\W^{1/2}\bSigma_{\beta}\W^{1/2}\Y_{i_{1}}\Y_{i_{1}}^{T}\W^{1/2}\bSigma_{\beta}\W^{1/2}\Y_{i})\right\}\\
&=2\sum_{\alpha=1}^{\beta-1}\sum_{i\neq i_{1}}^{m_{\alpha}}
\frac{d_{\alpha\beta}^{4}}{n_{\alpha}^{4}n_{\beta}^{4}}\rm{tr}(\W\bSigma_{\alpha}\W\bSigma_{\beta})^{2},
\end{align*}
(4)$\quad(\alpha=\alpha_{1}\neq\alpha_{2}=\alpha_{3}; i=i_{1}, i_{2}=i_{3})$
\begin{align*}
\E(I_{1}^{2})&=\E\bigg\{\sum_{\alpha\neq\alpha_{2}}^{\beta-1}\sum_{i=m_{\alpha-1}+1}^{m_{\alpha}}\sum_{i_{2}=m_{\alpha_{2}-1}+1}^{m_{\alpha_{2}}}
\frac{d_{\alpha\beta}^{2}d_{\alpha_{2}\beta}^{2}}{n_{\alpha}^{2}n_{\alpha_{2}}^{2}n_{\beta}^{4}}\\
&\quad(\Y_{i}^{T}\W^{1/2}\bSigma_{\beta}\W^{1/2}\Y_{i}\Y_{i_{2}}^{T}\W^{1/2}\bSigma_{\beta}\W^{1/2}\Y_{i_{2}})\bigg\}\\
&=\sum_{\alpha\neq\alpha_{2}}^{\beta-1}\sum_{i=m_{\alpha-1}+1}^{m_{\alpha}}\sum_{i_{2}=m_{\alpha_{2}-1}+1}^{m_{\alpha_{2}}}
\frac{d_{\alpha\beta}^{2}d_{\alpha_{2}\beta}^{2}}{n_{\alpha}^{2}n_{\alpha_{2}}^{2}n_{\beta}^{4}}\rm{tr}(\W\bSigma_{\alpha}\W\bSigma_{\beta})\rm{tr}(\W\bSigma_{\alpha_{2}}\W\bSigma_{\beta})\\
&=\sum_{\alpha\neq\alpha_{2}}^{\beta-1}\frac{d_{\alpha\beta}^{2}d_{\alpha_{2}\beta}^{2}}{n_{\alpha}n_{\alpha_{2}}n_{\beta}^{4}}\rm{tr}(\W\bSigma_{\alpha}\W\bSigma_{\beta})\rm{tr}(\W\bSigma_{\alpha_{2}}\W\bSigma_{\beta}),
\end{align*}
(5)$\quad(\alpha=\alpha_{2}\neq\alpha_{1}=\alpha_{3}; i=i_{2}, i_{1}=i_{3})$ or $(\alpha=\alpha_{3}\neq\alpha_{1}=\alpha_{2}; i=i_{3}, i_{1}=i_{2})$
\begin{align*}
\E(I_{1}^{2})&=2\E\bigg\{\sum_{\alpha\neq\alpha{1}}^{\beta-1}\sum_{i=m_{\alpha-1}+1}^{m_{\alpha}}\sum_{i_{1}=m_{\alpha_{1}-1}+1}^{m_{\alpha_{1}}}
\frac{d_{\alpha\beta}^{2}d_{\alpha_{1}\beta}^{2}}{n_{\alpha}^{2}n_{\alpha_{1}}^{2}n_{\beta}^{4}}\\
&\quad(\Y_{i}^{T}\W^{1/2}\bSigma_{\beta}\W^{1/2}\Y_{i_{1}}\Y_{i_{1}}^{T}\W^{1/2}\bSigma_{\beta}\W^{1/2}\Y_{i})\bigg\}\\
&=2\sum_{\alpha\neq\alpha{1}}^{\beta-1}\sum_{i=m_{\alpha-1}+1}^{m_{\alpha}}\sum_{i_{1}=m_{\alpha_{1}-1}+1}^{m_{\alpha_{1}}}
\frac{d_{\alpha\beta}^{2}d_{\alpha_{1}\beta}^{2}}{n_{\alpha}^{2}n_{\alpha_{1}}^{2}n_{\beta}^{4}}\rm{tr}(\W\bSigma_{\alpha}\W\bSigma_{\beta}\W\bSigma_{\alpha_{1}}\W\bSigma_{\beta})\\
&=2\sum_{\alpha\neq\alpha_{1}}^{\beta-1}\frac{d_{\alpha\beta}^{2}d_{\alpha_{1}\beta}^{2}}{n_{\alpha}n_{\alpha_{1}}n_{\beta}^{4}}\rm{tr}(\W\bSigma_{\alpha}\W\bSigma_{\beta}\W\bSigma_{\alpha_{1}}\W\bSigma_{\beta}).
\end{align*}
By (1), we have
\begin{align*}
I_{11}&=\sum_{\alpha=1}^{\beta-1}\sum_{i=m_{\alpha-1}+1}^{m_{\alpha}}
\frac{d_{\alpha\beta}^{4}\E(\Y_{i}^{T}\W^{1/2}\bSigma_{\beta}\W^{1/2}\Y_{i})^{2}}{n_{\alpha}^{4}n_{\beta}^{4}}.
\end{align*}
Combining (2) and (4), we get
 \begin{align*}
I_{121}&=\sum_{\alpha=1}^{\beta-1}\sum_{\alpha_{2}=1}^{\beta-1}\frac{d_{\alpha\beta}^{2}d_{\alpha_{2}\beta}^{2}\rm{tr}(\W\bSigma_{\alpha}\W\bSigma_{\beta})\rm{tr}(\W\bSigma_{\alpha_{2}}\W\bSigma_{\beta})}{n_{\alpha}n_{\alpha_{2}}n_{\beta}^{4}}
-\sum_{\alpha=1}^{\beta-1}\frac{d_{\alpha\beta}^{4}\rm{tr}^{2}(\W\bSigma_{\alpha}\W\bSigma_{\beta})}{n_{\alpha}^{3}n_{\beta}^{4}}.
\end{align*}
Merging (3) and (5), We have 
 \begin{align*}
I_{122}&=2\left\{\sum_{\alpha=1}^{\beta-1}\sum_{\alpha_{1}=1}^{\beta-1}\frac{d_{\alpha\beta}^{2}d_{\alpha_{1}\beta}^{2}\rm{tr}(\W\bSigma_{\alpha}\W\bSigma_{\beta}\W\bSigma_{\alpha_{1}}\W\bSigma_{\beta})}{n_{\alpha}n_{\alpha_{1}}n_{\beta}^{4}}
-\sum_{\alpha=1}^{\beta-1}\frac{d_{\alpha\beta}^{4}\rm{tr}(\W\bSigma_{\alpha}\W\bSigma_{\beta})^{2}}{n_{\alpha}^{3}n_{\beta}^{4}}\right\}.
\end{align*}
Using a splitting method similar to $\E(I_{1}^{2})$, $\E(I_{2}^{2})$ can be represented as $\E(I_{2}^{2})=I_{21}+I_{221}+I_{222}$, where
\begin{align*}
I_{21}&=\sum_{i=m_{\beta-1}+1}^{j-1}\frac{d_{\beta\beta}^{4}\E(\Y_{i}^{T}\W^{1/2}\bSigma_{\beta}\W^{1/2}\Y_{i})^{2}}{n_{\beta}^{4}(n_{\beta}-1)^{4}},
\end{align*}
 \begin{align*}
I_{221}&=\frac{(j-1-m_{\beta-1})^{2}}{n_{\beta}^{4}(n_{\beta}-1)^{4}}d_{\beta\beta}^{4}\rm{tr}^{2}(\W\bSigma_{\beta})^{2}
-\frac{(j-1-m_{\beta-1})}{n_{\beta}^{4}(n_{\beta}-1)^{4}}d_{\beta\beta}^{4}\rm{tr}^{2}(\W\bSigma_{\beta})^{2},
\end{align*}
and
\begin{align*}
I_{222}&=\sum_{i=m_{\beta-1}+1}^{j-1}\sum_{i\neq i{1}}\frac{2d_{\beta\beta}^{4}}{n_{\beta}^{4}(n_{\beta}-1)^{4}}\rm{tr}(\W\bSigma_{\beta})^{4}.
\end{align*}
Further, the rest are respectively 
\begin{align*}
\E(I_{3}^{2})&=\sum_{\alpha=1}^{\beta-1}\frac{(j-1-m_{\beta-1})}{n_{\alpha}n_{\beta}^{4}(n_{\beta}-1)^{2}}d_{\beta\beta}^{2}d_{\alpha\beta}^{2}
\rm{tr}(\W\bSigma_{\alpha}\W\bSigma_{\beta}\W\bSigma_{\beta}\W\bSigma_{\beta}),
\end{align*}
\begin{align*}
\E(I_{1}I_{2})&=\sum_{\alpha=1}^{\beta-1}\frac{(j-1-m_{\beta-1})}{n_{\alpha}n_{\beta}^{4}(n_{\beta}-1)^{2}}d_{\beta\beta}^{2}d_{\alpha\beta}^{2}
\rm{tr}(\W\bSigma_{\alpha}\W\bSigma_{\beta})\rm{tr}(\W\bSigma_{\beta})^{2},
\end{align*}
\begin{align*}
\E(I_{1}I_{1}^{'})&=\E(I_{1}^{2}),
\end{align*}
\begin{align*}
\E(I_{2}I_{2}^{'})&=\E(I_{2}^{2})+\sum_{i=m_{\beta-1}+1}^{j-1}\sum_{i_{2}=j}^{j_{1}-1}\frac{d_{\beta\beta}^{4}}{n_{\beta}^{4}(n_{\beta}-1)^{4}}\rm{tr}^{2}(\W\bSigma_{\beta})^{2},
\end{align*}
\begin{align*}
\E(I_{3}I_{3}^{'})&=\E(I_{3}I_{3}),
\end{align*}
\begin{align*}
\E(I_{1}I_{2}^{'})&=\sum_{\alpha=1}^{\beta-1}\frac{(j_{1}-1-m_{\beta-1})}{n_{\alpha}n_{\beta}^{4}(n_{\beta}-1)^{2}}d_{\beta\beta}^{2}d_{\alpha\beta}^{2}
\rm{tr}(\W\bSigma_{\alpha}\W\bSigma_{\beta})\rm{tr}(\W\bSigma_{\beta})^{2},
\end{align*}
\begin{align*}
\E(I_{2}I_{1}^{'})&=\sum_{\alpha=1}^{\beta-1}\frac{(j-1-m_{\beta-1})}{n_{\alpha}n_{\beta}^{4}(n_{\beta}-1)^{2}}d_{\beta\beta}^{2}d_{\alpha\beta}^{2}
\rm{tr}(\W\bSigma_{\alpha}\W\bSigma_{\beta})\rm{tr}(\W\bSigma_{\beta})^{2}.
\end{align*}
Furthermore, through \eqref{A.3}, we can obtain 
\begin{align*}
&\E^{2}\left(\sum_{j=m_{\beta-1}+1}^{m_{\beta}}\E(\V_{nj}^{2}|\mathcal{F}_{n,j-1})\right)\\
&=\left(\sum_{j=m_{\beta-1}+1}^{m_{\beta}}\E(\V_{nj}^{2})\right)^{2}\\
&=16\sum_{\alpha=1}^{\beta-1}\sum_{\alpha_{1}=1}^{\beta-1}\frac{d_{\alpha\beta}^{2}d_{\alpha_{1}\beta}^{2}}{n_{\alpha}n_{\alpha_{1}}n_{\beta}^{2}}
\rm{tr}(\W\bSigma_{\alpha}\W\bSigma_{\beta})\rm{tr}(\W\bSigma_{\alpha_{1}}\W\bSigma_{\beta})+4\frac{d_{\beta\beta}^{4}\rm{tr}^{2}(\W\bSigma_{\beta})^{2}}{n_{\beta}^{2}(n_{\beta}-1)^{2}}\\
&\quad+16\sum_{\alpha=1}^{\beta-1}\frac{d_{\alpha\beta}^{2}d_{\beta\beta}^{2}}{n_{\alpha}n_{\beta}^{2}(n_{\beta}-1)}\rm{tr}(\W\bSigma_{\alpha}\W\bSigma_{\beta})\rm{tr}(\W\bSigma_{\beta})^{2}.
\end{align*}
Therefore, by Assumption (A4), we can obtain 
\begin{align*}
&\Var\left(\sum_{j=m_{\beta-1}+1}^{m_{\beta}}\E(\V_{nj}^{2}|\mathcal{F}_{n,j-1})\right)\\
&=\E\left(\sum_{j=m_{\beta-1}+1}^{m_{\beta}}\E(\V_{nj}^{2}|\mathcal{F}_{n,j-1})\right)^{2}
-\E^{2}\left(\sum_{j=m_{\beta-1}+1}^{m_{\beta}}\E(\V_{nj}^{2}|\mathcal{F}_{n,j-1})\right)\\
&=\sum_{j=m_{\beta-1}+1}^{m_{\beta}}[16I_{11}+16I_{122}+16I_{21}+16I_{222}+64\E(I_{3}^{2})]\\
&\quad+2\sum_{j=m_{\beta-1}+1}^{m_{\beta}}\sum_{j<j_{1}}^{m_{\beta}}[16I_{11}+16I_{122}+16I_{21}+16I_{222}+64\E(I_{3}^{2})]
+o(\sigma^{4}({\rm{T}}_{n_0})).
\end{align*}
Further, by lemma \ref{lem 2.2}, we obtain 
\begin{align*}
&\sum_{j=m_{\beta-1}+1}^{m_{\beta}}16I_{11}+2\sum_{j=m_{\beta-1}+1}^{m_{\beta}}\sum_{j<j_{1}}^{m_{\beta}}16I_{11}\\
&=16\sum_{\alpha=1}^{\beta-1}\sum_{i=m_{\alpha-1}+1}^{m_{\alpha}}\frac{d_{\alpha\beta}^{4}}{n_{\alpha}^{4}n_{\beta}^{2}}
\E(\Y_{i}^{T}\W^{1/2}\bSigma_{\beta}\W^{1/2}\Y_{i})^{2}\\
&\leq 16\sum_{\alpha=1}^{\beta-1}\frac{d_{\alpha\beta}^{4}}{n_{\alpha}^{3}n_{\beta}^{2}}
\rm{tr}^{2}(\W\bSigma_{\alpha}\W\bSigma_{\beta})
+(2+\Delta)\sum_{\alpha=1}^{\beta-1}\frac{16d_{\alpha\beta}^{4}}{n_{\alpha}^{3}n_{\beta}^{2}}
\rm{tr}(\W\bSigma_{\alpha}\W\bSigma_{\beta})^{2},
\end{align*}
where the first term on the right is $o(\sigma^{4}({\rm{T}}_{n_0}))$, and  by Assumption (A4), the second term on the right is also $o(\sigma^{4}({\rm{T}}_{n_0}))$. By the same method, we can obtain 
\begin{align*}
&\sum_{j=m_{\beta-1}+1}^{m_{\beta}}16I_{122}+2\sum_{j=m_{\beta-1}+1}^{m_{\beta}}\sum_{j<j_{1}}^{m_{\beta}}16I_{122}=o(\sigma^{4}({\rm{T}}_{n_0})),\\
&\sum_{j=m_{\beta-1}+1}^{m_{\beta}}16I_{21}+2\sum_{j=m_{\beta-1}+1}^{m_{\beta}}\sum_{j<j_{1}}^{m_{\beta}}16I_{21}=o(\sigma^{4}({\rm{T}}_{n_0})),\\
&\sum_{j=m_{\beta-1}+1}^{m_{\beta}}16I_{222}+2\sum_{j=m_{\beta-1}+1}^{m_{\beta}}\sum_{j<j_{1}}^{m_{\beta}}16I_{222}=o(\sigma^{4}({\rm{T}}_{n_0})),\\
&\sum_{j=m_{\beta-1}+1}^{m_{\beta}}64\E(I_{3}^{2})+2\sum_{j=m_{\beta-1}+1}^{m_{\beta}}\sum_{j<j_{1}}^{m_{\beta}}64\E(I_{3}^{2})=o(\sigma^{4}({\rm{T}}_{n_0})).
\end{align*}
Therefore, we have 
\begin{align}\label{A.4}
\Var\left(\sum_{j=m_{\alpha-1}+1}^{m_{\alpha}}\E(\V_{nj}^{2}|\mathcal{F}_{n,j-1})\right)=o(\sigma^{4}({\rm{T}}_{n_0})),
\end{align}
for $\alpha=1,\ldots,K.$

For proving $\sum_{j=1}^{n}\E(\V_{nj}^{2}|\mathcal{F}_{n,j-1})/\sigma^{2}({\rm{T}}_{n_0})\stackrel{p}{\longrightarrow}1$, we only need to prove 
$$\Var\left(\sum_{j=1}^{n}\E(\V_{nj}^{2}|\mathcal{F}_{n,j-1})\right)=o(\sigma^{4}({\rm{T}}_{n_0})),$$
where 
\begin{align*}
\sigma^{4}({\rm{T}}_{n_0})&=\E^{2}\left(\sum_{j=1}^{n}\E(\V_{nj}^{2}|\mathcal{F}_{n,j-1})\right)\\
&=16\sum_{\beta=1}^{K}\sum_{\alpha=1}^{\beta-1}\sum_{\beta_{1}=1}^{K}\sum_{\alpha_{1}=1}^{\beta_{1}-1}
\frac{d_{\alpha\beta}^{2}d_{\alpha_{1}\beta_{1}}^{2}}{n_{\alpha}n_{\beta}n_{\alpha_{1}}n_{\beta_{1}}}
\rm{tr}(\W\bSigma_{\alpha}\W\bSigma_{\beta})\rm{tr}(\W\bSigma_{\alpha_{1}}\W\bSigma_{\beta_{1}})\\
&\quad+4\sum_{\beta=1}^{K}\sum_{\beta_{1}=1}^{K}\frac{d_{\beta\beta}^{2}d_{\beta_{1}\beta_{1}}^{2}\rm{tr}(\W\bSigma_{\beta})^{2}
\rm{tr}(\W\bSigma_{\beta_{1}})^{2}}{n_{\beta}(n_{\beta}-1)n_{\beta_{1}}(n_{\beta_{1}}-1)}\\
&\quad+16\sum_{\beta=1}^{K}\sum_{\alpha=1}^{\beta-1}\frac{d_{\alpha\beta}^{2}}{n_{\alpha}n_{\beta}}\rm{tr}(\W\bSigma_{\alpha}\W\bSigma_{\beta})
\sum_{\beta_{1}=1}^{K}\frac{d_{\beta_{1}\beta_{1}}^{2}\rm{tr}(\W\bSigma_{\beta_{1}})^{2}}{n_{\beta_{1}}(n_{\beta_{1}}-1)}.
\end{align*}
Since 
\begin{align*}
\Var\left(\sum_{j=1}^{n}\E(\V_{nj}^{2}|\mathcal{F}_{n,j-1})\right)
&=\Var\left(\sum_{\beta=1}^{K}\sum_{j=m_{\beta-1}+1}^{m_{\beta}}\E(\V_{nj}^{2}|\mathcal{F}_{n,j-1})\right)\\
&\leq K\sum_{\beta=1}^{K}\Var\left(\sum_{j=m_{\beta-1}+1}^{m_{\beta}}\E(\V_{nj}^{2}|\mathcal{F}_{n,j-1})\right),
\end{align*}
we get $\Var\left(\sum_{j=1}^{n}\E(\V_{nj}^{2}|\mathcal{F}_{n,j-1})\right)=o(\sigma^{4}({\rm{T}}_{n_0}))$ from \eqref{A.4}. Therefore, we have $$\sum_{j=1}^{n}\E(\V_{nj}^{2}|\mathcal{F}_{n,j-1})/\sigma^{2}({\rm{T}}_{n_0})\stackrel{p}{\longrightarrow}1.$$

Next, we will prove that $\sum_{j=1}^{n}\E(\V_{nj}^{4})=o(\sigma^{4}({\rm{T}}_{n_0})).$ By \eqref{A.1}, we can obtain 
\begin{align}\label{A.5}
\sum_{j=1}^{n}\E(\V_{nj}^{4})&=\sum_{\beta=1}^{K}\sum_{j=m_{\beta-1}+1}^{m_{\beta}}\E[(\V_{nj}^{2})^{2}]\nonumber\\
&=\sum_{\beta=1}^{K}\sum_{j=m_{\beta-1}+1}^{m_{\beta}}\E\bigg[4\left(\sum_{\alpha=1}^{\beta-1}
\sum_{i=m_{\alpha-1}+1}^{m_{\alpha}}\sum_{\alpha_{1}=1}^{\beta-1}\sum_{i_{1}=m_{\alpha_{1}-1}+1}^{m_{\alpha_{1}}}
\frac{d_{\alpha\beta}d_{\alpha_{1}\beta}}{n_{\alpha}n_{\alpha_{1}}n_{\beta}^{2}}\Y_{j}^{T}\Y_{i}\Y_{i_{1}}^{T}\Y_{j}\right)\nonumber\\
&\quad+4\left(\sum_{i=m_{\beta-1}+1}^{j-1}\sum_{i_{1}=m_{\beta-1}+1}^{j-1}\frac{d_{\beta\beta}^{2}}{n_{\beta}^{2}(n_{\beta}-1)^{2}}\Y_{j}^{T}\Y_{i}\Y_{i_{1}}^{T}\Y_{j}\right)\nonumber\\
&\quad+8\left(\sum_{\alpha=1}^{\beta-1}\sum_{i=m_{\alpha-1}+1}^{m_{\alpha}}\sum_{i_{1}=m_{\beta-1}+1}^{j-1}
\frac{d_{\alpha\beta}d_{\beta\beta}}{n_{\alpha}n_{\beta}^{2}(n_{\beta}-1)}\Y_{j}^{T}\Y_{i}\Y_{i_{1}}^{T}\Y_{j}\right)\bigg]^{2}\nonumber\\
&=16\sum_{\beta=1}^{K}\sum_{j=m_{\beta-1}+1}^{m_{\beta}}\sum_{\alpha=1}^{\beta-1}
\frac{d_{\alpha\beta}^{4}\E(\Y_{i}^{T}\Y_{j})^{4}}{n_{\alpha}^{3}n_{\beta}^{4}}\nonumber\\
&\quad+48\sum_{\beta=1}^{K}\sum_{j=m_{\beta-1}+1}^{m_{\beta}}\sum_{\alpha=1}^{\beta-1}\sum_{i=m_{\alpha-1}+1}^{m_{\alpha}}\sum_{\alpha_{1}=1}^{\beta-1}
\sum_{i\neq i_{1}}^{m_{\alpha_{1}}}\frac{d_{\alpha\beta}^{2}d_{\alpha_{1}\beta}^{2}\E\left[(\Y_{i}^{T}\Y_{j})^{2}
(\Y_{i_{1}}^{T}\Y_{j})^{2}\right]}{n_{\alpha}^{2}n_{\alpha_{1}}^{2}n_{\beta}^{4}}\nonumber\\
&\quad+96\sum_{\beta=1}^{K}\sum_{j=m_{\beta-1}+1}^{m_{\beta}}\sum_{\alpha=1}^{\beta-1}\sum_{i=m_{\alpha-1}+1}^{m_{\alpha}}\sum_{i_{1}=m_{\beta-1}+1}^{j-1}
\frac{d_{\alpha\beta}^{2}d_{\beta\beta}^{2}\E\left[(\Y_{i}^{T}\Y_{j})^{2}
(\Y_{i_{1}}^{T}\Y_{j})^{2}\right]}{n_{\alpha}^{2}n_{\beta}^{4}(n_{\beta}-1)^{2}}\nonumber\\
&\quad+16\sum_{\beta=1}^{K}\sum_{j=m_{\beta-1}+1}^{m_{\beta}}\sum_{i=m_{\beta-1}+1}^{j-1}
\frac{d_{\beta\beta}^{4}\E(\Y_{i}^{T}\Y_{j})^{4}}{n_{\beta}^{4}(n_{\beta}-1)^{4}}\nonumber\\
&\quad+48\sum_{\beta=1}^{K}\sum_{j=m_{\beta-1}+1}^{m_{\beta}}\sum_{i=m_{\beta-1}+1}^{j-1}\sum_{i\neq i_{1}}
\frac{d_{\beta\beta}^{4}\E\left[(\Y_{i}^{T}\Y_{j})^{2}(\Y_{i_{1}}^{T}\Y_{j})^{2}\right]}{n_{\beta}^{4}(n_{\beta}-1)^{4}}.
\end{align}

Since 
\begin{align*}
&\sum_{j=m_{\beta-1}+1}^{m_{\beta}}\sum_{\alpha=1}^{\beta-1}\sum_{i=m_{\alpha-1}+1}^{m_{\alpha}}\sum_{\alpha_{1}=1}^{\beta-1}
\sum_{i\neq i_{1}}^{m_{\alpha_{1}}}\frac{d_{\alpha\beta}^{2}d_{\alpha_{1}\beta}^{2}\E\left[(\Y_{i}^{T}\Y_{j})^{2}
(\Y_{i_{1}}^{T}\Y_{j})^{2}\right]}{n_{\alpha}^{2}n_{\alpha_{1}}^{2}n_{\beta}^{4}}\\
&\quad=\sum_{j=m_{\beta-1}+1}^{m_{\beta}}\sum_{\alpha=1}^{\beta-1}\sum_{i=m_{\alpha-1}+1}^{m_{\alpha}}\sum_{\alpha_{1}=1}^{\beta-1}
\sum_{i_{1}=m_{\alpha_{1}-1}}^{m_{\alpha_{1}}}\frac{d_{\alpha\beta}^{2}d_{\alpha_{1}\beta}^{2}\E\left[(\Y_{i}^{T}\Y_{j})^{2}
(\Y_{i_{1}}^{T}\Y_{j})^{2}\right]}{n_{\alpha}^{2}n_{\alpha_{1}}^{2}n_{\beta}^{4}}\\
&\quad-\sum_{j=m_{\beta-1}+1}^{m_{\beta}}\sum_{\alpha=1}^{\beta-1}\sum_{i=m_{\alpha-1}+1}^{m_{\alpha}}
\frac{d_{\alpha\beta}^{4}E\left[(\Y_{i}^{T}\Y_{j})^{4}\right]}{n_{\alpha}^{4}n_{\beta}^{4}},
\end{align*}
by lemma \ref{lem 2.2} and Assumption (A4), the first term on the right side of the above equation is $o(\sigma^{4}({\rm{T}}_{n_0}))$, and by lemma \ref{lem 2.3} and Assumption (A4), the second term is also $o(\sigma^{4}({\rm{T}}_{n_0}))$. Similarly, we can obtain 

\begin{align*}
&\sum_{\beta=1}^{K}\sum_{j=m_{\beta-1}+1}^{m_{\beta}}\sum_{\alpha=1}^{\beta-1}\frac{d_{\alpha\beta}^{4}\E(\Y_{i}^{T}\Y_{j})^{4}}{n_{\alpha}^{3}n_{\beta}^{4}}
=o(\sigma^{4}({\rm{T}}_{n_0})),\\
&\sum_{\beta=1}^{K}\sum_{j=m_{\beta-1}+1}^{m_{\beta}}\sum_{\alpha=1}^{\beta-1}\sum_{i=m_{\alpha-1}+1}^{m_{\alpha}}\sum_{i_{1}=m_{\beta-1}+1}^{j-1}
\frac{d_{\alpha\beta}^{2}d_{\beta\beta}^{2}\E\left[(\Y_{i}^{T}\Y_{j})^{2}
(\Y_{i_{1}}^{T}\Y_{j})^{2}\right]}{n_{\alpha}^{2}n_{\beta}^{4}(n_{\beta}-1)^{2}}=o(\sigma^{4}({\rm{T}}_{n_0})),\\
&\sum_{\beta=1}^{K}\sum_{j=m_{\beta-1}+1}^{m_{\beta}}\sum_{i=m_{\beta-1}+1}^{j-1}
\frac{d_{\beta\beta}^{4}\E(\Y_{i}^{T}\Y_{j})^{4}}{n_{\beta}^{4}(n_{\beta}-1)^{4}}=o(\sigma^{4}({\rm{T}}_{n_0})),\\
&\sum_{\beta=1}^{K}\sum_{j=m_{\beta-1}+1}^{m_{\beta}}\sum_{i=m_{\beta-1}+1}^{j-1}\sum_{i\neq i_{1}}
\frac{d_{\beta\beta}^{4}\E\left[(\Y_{i}^{T}\Y_{j})^{2}(\Y_{i_{1}}^{T}\Y_{j})^{2}\right]}{n_{\beta}^{4}(n_{\beta}-1)^{4}}=o(\sigma^{4}({\rm{T}}_{n_0})).
\end{align*}
Combining these with \eqref{A.5}, we have $\sum_{j=1}^{n}\E(\V_{nj}^{4})=o(\sigma^{4}({\rm{T}}_{n_0}))$. Then the lemma proof is complete.

\subsection*{Proof of theorem 2.2}
\noindent Note that ${\rm{T}}_{n_0}=\sum_{j=1}^{n}\V_{nj}$. By Lemma \ref{lem 2.4}, we have 
$$\sigma^{-4}({\rm{T}}_{n_0})\sum_{j=1}^{n}\E(\V_{nj}^{4})\rightarrow 0,$$
therefore we can obtain 
\begin{align}\label{A.6}
\sigma^{-4}({\rm{T}}_{n_0})\sum_{j=1}^{n}\E(\V_{nj}^{4}|\mathcal{F}_{n,j-1})\stackrel{p}{\longrightarrow} 0.
\end{align}
Further, by Markov's inequality 
$$\sum_{j=1}^{n}\sigma^{-2}({\rm{T}}_{n_0})\E(\V_{nj}^{2}\I(|\V_{nj}|>\varepsilon\sigma({\rm{T}}_{n_0}))|\mathcal{F}_{n,j-1})\leq\sigma^{-4}({\rm{T}}_{n_0})
\varepsilon^{-2}\sum_{j=1}^{n}\E(\V_{nj}^{4}|\mathcal{F}_{n,j-1}), \forall\varepsilon>0$$ 
and \eqref{A.6}, we can find that Lindbergh's condition holds 
$$\sum_{j=1}^{n}\sigma^{-2}({\rm{T}}_{n_0})\E(\V_{nj}^{2}\I(|\V_{nj}|>\varepsilon\sigma({\rm{T}}_{n_0}))|\mathcal{F}_{n,j-1})\stackrel{p}{\longrightarrow}0.$$ 
Again using lemma \ref{lem 2.4}, we have
$$\frac{\sum_{j=1}^{n}\E(\V_{nj}^{2}|\mathcal{F}_{n,j-1})}{\sigma^{2}({\rm{T}}_{n_0})}\stackrel{p}{\longrightarrow}1.$$ 
So, using the martingale difference central limit theorem, we have
$$\frac{{\rm{T}}_{n_0}}{\sigma({\rm{T}}_{n_0})}\stackrel{d}{\longrightarrow}N(0,1).$$
Then the theorem is proved.

\subsection*{Proof of theorem 2.3}
\noindent Firstly, note that 
\begin{align*}
\overline{\X}_{\alpha} &=\frac{1}{n_\alpha}\sum_{i=1}^{n_\alpha}\X_{\alpha i}=\overline{\Y}_{\alpha}-\bmu_{\alpha},\\
\widehat{\bSigma_{\alpha}}&=\frac{1}{n_\alpha-1}\sum_{i=1}^{n_\alpha}(\X_{\alpha i}-\overline{\X}_{\alpha})(\X_{\alpha i}-\overline{\X}_{\alpha})^T\\ 
&= \frac{1}{n_\alpha-1}\sum_{i=1}^{n_\alpha}(\Y_{\alpha i}-\overline{\Y}_{\alpha})(\Y_{\alpha i}-\overline{\Y}_{\alpha})^T.                                            
\end{align*}
So, we have
\begin{align*}
\W\widehat{\bSigma_{\alpha}}&=\frac{1}{n_\alpha-1}\sum_{i=1}^{n_\alpha}\W(\Y_{\alpha i}-\overline{\Y}_{\alpha})(\Y_{\alpha i}-\overline{\Y}_{\alpha})^T\\
&=\frac{1}{n_\alpha-1}\sum_{i=1}^{n_\alpha}\W\left(\frac{n_{\alpha}-1}{n_{\alpha}}\Y_{\alpha i}-
\frac{1}{n_{\alpha}}\sum_{\substack{j=1\\ i\neq j}}^{n_{\alpha}}{\Y}_{\alpha j}\right)\left(\frac{n_{\alpha}-1}{n_{\alpha}}\Y_{\alpha i}-
\frac{1}{n_{\alpha}}\sum_{\substack{j=1\\ i\neq j}}^{n_{\alpha}}{\Y}_{\alpha j}\right)^T\\
&=\frac{1}{n_\alpha}\sum_{i=1}^{n_\alpha}\W\Y_{\alpha i}\Y_{\alpha i}^T-\frac{1}{n_{\alpha}(n_{\alpha}-1)}\sum_{\substack{i,j=1\\ i\neq j}}^{n_{\alpha}}
\W\Y_{\alpha i}\Y_{\alpha j}^T.
\end{align*}
From the above equation, we can have 
\begin{align*}
{\rm{tr}}(\W\widehat{\bSigma_{\alpha}})^{2}&=\frac{1}{{n_{\alpha}}^2}\sum_{i=1}^{n_\alpha}(\Y_{\alpha i}^T\W\Y_{\alpha i})^{2}+
\frac{1}{{n_{\alpha}}^2}\sum_{\substack{i,j=1\\ i\neq j}}^{n_{\alpha}}(\Y_{\alpha i}^T\W\Y_{\alpha j})^{2}\\
&\quad+\frac{1}{n_{\alpha}^{2}(n_{\alpha}-1)^{2}}\sum_{\substack{i,j=1\\ i\neq j}}^{n_{\alpha}}\sum_{\substack{k,l=1\\ k\neq l}}^{n_{\alpha}}
\Y_{\alpha i}^T\W\Y_{\alpha l}\Y_{\alpha k}^T\W\Y_{\alpha j}\\
&\quad-\frac{1}{n_{\alpha}^{2}(n_{\alpha}-1)}\sum_{\substack{i,j,k=1\\ k\neq j}}^{n_{\alpha}}
\Y_{\alpha j}^T\W\Y_{\alpha i}\Y_{\alpha i}^T\W\Y_{\alpha k}\\
&=\frac{1}{{n_{\alpha}}^2}\sum_{i=1}^{n_\alpha}(\Y_{\alpha i}^T\W\Y_{\alpha i})^{2}+\frac{n_{\alpha}^{2}-2n_{\alpha}+2}{n_{\alpha}^{2}(n_{\alpha}-1)^{2}}
\sum_{\substack{i,j=1\\ i\neq j}}^{n_{\alpha}}(\Y_{\alpha i}^T\W\Y_{\alpha j})^{2}\\
&\quad+\frac{1}{n_{\alpha}^{2}(n_{\alpha}-1)^{2}}\sum_{\substack{i,j=1\\ i\neq j}}^{n_{\alpha}}\Y_{\alpha i}^T\W\Y_{\alpha i}\Y_{\alpha j}^T\W\Y_{\alpha j}\\
&\quad-\frac{4}{n_{\alpha}^{2}(n_{\alpha}-1)}\sum_{\substack{i,j=1\\ i\neq j}}^{n_{\alpha}}\Y_{\alpha i}^T\W\Y_{\alpha i}\Y_{\alpha i}^T\W\Y_{\alpha j}\\
&\quad+\frac{2}{n_{\alpha}^{2}(n_{\alpha}-1)^{2}}\sum_{\substack{i,j,k=1\\ i\neq j,k\neq j,k\neq i}}^{n_{\alpha}}\Y_{\alpha i}^T\W\Y_{\alpha i}\Y_{\alpha j}^T\W\Y_{\alpha k}\\
&\quad-\frac{2n_{\alpha}-4}{n_{\alpha}^{2}(n_{\alpha}-1)^{2}}\sum_{\substack{i,j,k=1\\ i\neq j,k\neq j,k\neq i}}^{n_{\alpha}}\Y_{\alpha i}^T\W\Y_{\alpha j}\Y_{\alpha i}^T\W\Y_{\alpha k}\\
&\quad+\frac{1}{n_{\alpha}^{2}(n_{\alpha}-1)^{2}}\sum_{\substack{i,j,k,l=1\\ i\neq j\neq k\neq l}}^{n_{\alpha}}\Y_{\alpha i}^T\W\Y_{\alpha l}\Y_{\alpha k}^T\W\Y_{\alpha j},
\end{align*}

\begin{align*}
[{\rm{tr}}(\W\widehat{\bSigma_{\alpha}})]^{2}&=\bigg[\frac{1}{n_\alpha}\sum_{i=1}^{n_\alpha}\Y_{\alpha i}^T\W\Y_{\alpha i}-\frac{1}{n_{\alpha}(n_{\alpha}-1)}\sum_{\substack{i,j=1\\ i\neq j}}^{n_{\alpha}}\Y_{\alpha i}^T\W\Y_{\alpha j}\bigg]^2\\
&=\frac{1}{{n_{\alpha}}^2}\sum_{i=1}^{n_\alpha}(\Y_{\alpha i}^T\W\Y_{\alpha i})^{2}+\frac{1}{{n_{\alpha}}^2}\sum_{\substack{i,j=1\\ i\neq j}}^{n_{\alpha}}\Y_{\alpha i}^T\W\Y_{\alpha i}\Y_{\alpha j}^T\W\Y_{\alpha j}\\
&\quad+\frac{2}{n_{\alpha}^{2}(n_{\alpha}-1)^{2}}\sum_{\substack{i,j=1\\ i\neq j}}^{n_{\alpha}}(\Y_{\alpha i}^T\W\Y_{\alpha j})^{2}\\
&\quad-\frac{4}{n_{\alpha}^{2}(n_{\alpha}-1)}\sum_{\substack{i,j=1\\ i\neq j}}^{n_{\alpha}}\Y_{\alpha i}^T\W\Y_{\alpha i}\Y_{\alpha i}^T\W\Y_{\alpha j}\\
&\quad-\frac{2}{n_{\alpha}^{2}(n_{\alpha}-1)}\sum_{\substack{i,j,k=1\\ i\neq j,k\neq j,k\neq i}}^{n_{\alpha}}\Y_{\alpha i}^T\W\Y_{\alpha i}\Y_{\alpha j}^T\W\Y_{\alpha k}\\
&\quad+\frac{4}{n_{\alpha}^{2}(n_{\alpha}-1)^{2}}\sum_{\substack{i,j,k=1\\ i\neq j,k\neq j,k\neq i}}^{n_{\alpha}}\Y_{\alpha i}^T\W\Y_{\alpha j}\Y_{\alpha i}^T\W\Y_{\alpha k}\\
&\quad+\frac{1}{n_{\alpha}^{2}(n_{\alpha}-1)^{2}}\sum_{\substack{i,j,k,l=1\\ i\neq j\neq k\neq l}}^{n_{\alpha}}\Y_{\alpha i}^T\W\Y_{\alpha l}\Y_{\alpha k}^T\W\Y_{\alpha j},
\end{align*}
and
\begin{align*}
\Q_{\alpha}&=\frac{1}{(n_{\alpha}-1)}\sum_{i=1}^{n_{\alpha}}{\lVert\W^{1/2}(\Y_{\alpha i}-\overline{\Y}_{\alpha})\rVert}^4\\
&=\frac{1}{(n_{\alpha}-1)}\sum_{i=1}^{n_{\alpha}}\bigg[(\Y_{\alpha i}^T\W\Y_{\alpha i})^{2}+4(\Y_{\alpha i}^T\W\overline{\Y}_{\alpha})^{2}
+(\overline{\Y}_{\alpha}^{T}\W\overline{\Y}_{\alpha})^{2}\\
&\quad+2\Y_{\alpha i}^T\W\Y_{\alpha i}\overline{\Y}_{\alpha}^{T}\W\overline{\Y}_{\alpha}-4\Y_{\alpha i}^T\W\Y_{\alpha i}\Y_{\alpha i}^T\W\overline{\Y}_{\alpha}-4\Y_{\alpha i}^T\W\overline{\Y}_{\alpha}\overline{\Y}_{\alpha}^{T}\W\overline{\Y}_{\alpha}\bigg]\\
&=\frac{n_{\alpha}^{2}-3n_{\alpha}+3}{n_{\alpha}^{3}}\sum_{i=1}^{n_\alpha}(\Y_{\alpha i}^T\W\Y_{\alpha i})^{2}+
\frac{4n_{\alpha}-6}{n_{\alpha}^{3}(n_{\alpha}-1)}\sum_{\substack{i,j=1\\ i\neq j}}^{n_{\alpha}}(\Y_{\alpha i}^T\W\Y_{\alpha j})^{2}\\
&\quad+\frac{2n_{\alpha}-3}{n_{\alpha}^{3}(n_{\alpha}-1)}\sum_{\substack{i,j=1\\ i\neq j}}^{n_{\alpha}}\Y_{\alpha i}^T\W\Y_{\alpha i}\Y_{\alpha j}^T\W\Y_{\alpha j}\\
&\quad-\frac{4(n_{\alpha}^{2}-3n_{\alpha}+3)}{n_{\alpha}^{3}(n_{\alpha}-1)}\sum_{\substack{i,j=1\\ i\neq j}}^{n_{\alpha}}\Y_{\alpha i}^T\W\Y_{\alpha i}\Y_{\alpha i}^T\W\Y_{\alpha j}\\
&\quad+\frac{2n_{\alpha}-6}{n_{\alpha}^{3}(n_{\alpha}-1)}\sum_{\substack{i,j,k=1\\ i\neq j,k\neq j,k\neq i}}^{n_{\alpha}}\Y_{\alpha i}^T\W\Y_{\alpha i}\Y_{\alpha j}^T\W\Y_{\alpha k}\\
&\quad+\frac{4n_{\alpha}-12}{n_{\alpha}^{3}(n_{\alpha}-1)}\sum_{\substack{i,j,k=1\\ i\neq j,k\neq j,k\neq i}}^{n_{\alpha}}\Y_{\alpha i}^T\W\Y_{\alpha j}\Y_{\alpha i}^T\W\Y_{\alpha k}\\
&\quad-\frac{3}{n_{\alpha}^{3}(n_{\alpha}-1)}\sum_{\substack{i,j,k,l=1\\ i\neq j\neq k\neq l}}^{n_{\alpha}}\Y_{\alpha i}^T\W\Y_{\alpha l}\Y_{\alpha k}^T\W\Y_{\alpha j}.
\end{align*}
Therefore, we have 
\begin{align}\label{A.7}
{\rm{tr}}\widehat{(\W\bSigma_{\alpha})^2}&=\frac{n_{\alpha}-1}{n_{\alpha}(n_{\alpha}-2)(n_{\alpha}-3)}\left((n_{\alpha}-1)(n_{\alpha}-2){\rm{tr}}(\W\widehat\bSigma_{\alpha})^{2
}+{\rm{tr}}^{2}(\W\widehat\bSigma_{\alpha})-n_{\alpha}\Q_{\alpha}\right)\nonumber\\
&=\sum_{ i\neq j}^{n_{\alpha}}\frac{(\Y_{\alpha i}^T\W\Y_{\alpha j})^{2}}{n_{\alpha}(n_{\alpha}-1)}-2\sum_{i\neq j,k\neq j,k\neq i}^{n_{\alpha}}
\frac{\Y_{\alpha i}^T\W\Y_{\alpha j}\Y_{\alpha j}^T\W\Y_{\alpha k}}{n_{\alpha}(n_{\alpha}-1)(n_{\alpha}-2)}\nonumber\\
&\quad+\sum_{ i\neq j\neq k\neq l}^{n_{\alpha}}\frac{\Y_{\alpha i}^T\W\Y_{\alpha j}\Y_{\alpha k}^T\W\Y_{\alpha l}}{n_{\alpha}(n_{\alpha}-1)(n_{\alpha}-2)(n_{\alpha}-3)}\nonumber\\
&:=\Y_{2,n}^{\alpha}-2\Y_{4,n}^{\alpha}+\Y_{5,n}^{\alpha}.
\end{align}
From this, we have $\E\left({\rm{tr}}\widehat{(\W\bSigma_{\alpha})^2}\right)={\rm{tr}}(\W\bSigma_{\alpha})^2$, and now we only need to prove that 
$\Var\left({\rm{tr}}\widehat{(\W\bSigma_{\alpha})^2}\right)=o({\rm{tr}}^{2}(\W\bSigma_{\alpha})^2)$.
Further, we can conclude that $\Var\left({\rm{tr}}\widehat{(\W\bSigma_{\alpha})^2}\right)\leq3\left(\Var(\Y_{2,n}^{\alpha})+4\Var({\Y_{4,n}^{\alpha}})+\Var(\Y_{5,n}^{\alpha})\right)$. Then under given conditions, apply Proposition A.2. of Chen et al.(2010) when in our context, we have 
\begin{align*}
\Var(\Y_{2,n}^{\alpha})&=\frac{4{\rm{tr}}^{2}(\W\bSigma_{\alpha})^2}{n_{\alpha}^{2}}+\frac{8{\rm{tr}}(\W\bSigma_{\alpha})^4}{n_{\alpha}}+
\frac{4\Delta{\rm{tr}}(\A_{\alpha}^2\circ\A_{\alpha}^2)}{n_{\alpha}}\\
&\quad+O(\frac{{\rm{tr}}^{2}(\W\bSigma_{\alpha})^2}{n_{\alpha}^{3}}
+\frac{{\rm{tr}}(\W\bSigma_{\alpha})^4}{n_{\alpha}^{2}}),\\
\Var({\Y_{4,n}^{\alpha}})&=\frac{2{\rm{tr}}^{2}(\W\bSigma_{\alpha})^2}{n_{\alpha}^{3}}+\frac{2{\rm{tr}}(\W\bSigma_{\alpha})^4}{n_{\alpha}^2}+O(\frac{{\rm{tr}}^{2}(\W\bSigma_{\alpha})^2}{n_{\alpha}^{4}}
+\frac{{\rm{tr}}(\W\bSigma_{\alpha})^4}{n_{\alpha}^{3}}),\\
\Var({\Y_{5,n}^{\alpha}})&=\frac{8{\rm{tr}}^{2}(\W\bSigma_{\alpha})^2}{n_{\alpha}^{4}}+O(\frac{{\rm{tr}}^{2}(\W\bSigma_{\alpha})^2}{n_{\alpha}^{5}}
+\frac{{\rm{tr}}(\W\bSigma_{\alpha})^4}{n_{\alpha}^{4}}),
\end{align*}
where $\A_{\alpha}=\bGam_{\alpha}^T\W\bGam_{\alpha}.$ Since ${\rm{tr}}(\A_{\alpha}^2\circ\A_{\alpha}^2)\leq{\rm{tr}}(\W\bSigma_{\alpha})^4 $, and ${\rm{tr}}(\W\bSigma_{\alpha})^4=o({\rm{tr}}^{2}(\W\bSigma_{\alpha})^2)$ by Assumption (A4), we have $\Var\left({\rm{tr}}\widehat{(\W\bSigma_{\alpha})^2}\right)=o({\rm{tr}}^{2}(\W\bSigma_{\alpha})^2)$. So, $\frac{{\rm{tr}}\widehat{(\W\bSigma_{\alpha})^2}}{{\rm{tr}}(\W\bSigma_{\alpha})^{2}}\stackrel{p}{\longrightarrow}1$.

Further, we have 
\begin{align*}
&{\rm{tr}}(\W\widehat{\bSigma_{\alpha}}\W\widehat{\bSigma_{\beta}})\\
&={\rm{tr}}\bigg\{\frac{1}{(n_\alpha-1)(n_\beta-1)}\sum_{i=1}^{n_\alpha}\W(\Y_{\alpha i}-\overline{\Y}_{\alpha})(\Y_{\alpha i}-\overline{\Y}_{\alpha})^T\\
&\quad \sum_{j=1}^{n_{\beta}}\W(\Y_{{\beta} j}-\overline{\Y}_{\beta})(\Y_{\beta j}-\overline{\Y}_{\beta})^T\bigg\}\\
&=\frac{1}{n_\alpha n_\beta}{\rm{tr}}\bigg\{\sum_{i=1}^{n_\alpha}\sum_{j=1}^{n_{\beta}}\W^{1/2}(\Y_{\alpha i}-\overline{\Y}_{\alpha(i)})\Y_{\alpha i}^T\W
(\Y_{{\beta} j}-\overline{\Y}_{\beta(j)})\Y_{{\beta} j}^T\W^{1/2}\bigg\}\\
&=\frac{1}{n_\alpha n_\beta}{\rm{tr}}\bigg\{\sum_{i=1}^{n_\alpha}\sum_{j=1}^{n_{\beta}}\W^{1/2}(\Y_{\alpha i}-\bmu_{\alpha}+\bmu_{\alpha}-\overline{\Y}_{\alpha(i)})(\Y_{\alpha i}-\bmu_{\alpha}+\bmu_{\alpha})^T\\
&\quad\W(\Y_{{\beta} j}-\bmu_{\beta}+\bmu_{\beta}-\overline{\Y}_{\beta(j)})(\Y_{{\beta} j}-\bmu_{\beta}+\bmu_{\beta})^T\W^{1/2}\bigg\}\\
&:=\sum_{l=1}^{16}D_{l},
\end{align*}
where $\overline{\Y}_{\beta(j)}$ represents the $\beta$th sample mean of after excluding $\Y_{\beta j}$, and
\begin{align*}
&D_{1}=\frac{1}{n_\alpha n_\beta}\sum_{i=1}^{n_\alpha}\sum_{j=1}^{n_{\beta}}{\rm{tr}}\bigg\{\W^{1/2}(\Y_{\alpha i}-\bmu_{\alpha})(\Y_{\alpha i}-\bmu_{\alpha})^T\W(\Y_{{\beta} j}-\bmu_{\beta})(\Y_{{\beta} j}-\bmu_{\beta})^T\W^{1/2}\bigg\},\\
&D_{2}=\frac{-1}{n_\alpha n_\beta}\sum_{i=1}^{n_\alpha}\sum_{j=1}^{n_{\beta}}{\rm{tr}}\bigg\{\W^{1/2}(\overline{\Y}_{\alpha(i)}-\bmu_{\alpha})(\Y_{\alpha i}-\bmu_{\alpha})^T\W(\Y_{{\beta} j}-\bmu_{\beta})(\Y_{{\beta} j}-\bmu_{\beta})^T\W^{1/2}\bigg\},\\
&D_{3}=\frac{-1}{n_\alpha n_\beta}\sum_{i=1}^{n_\alpha}\sum_{j=1}^{n_{\beta}}{\rm{tr}}\bigg\{\W^{1/2}(\Y_{\alpha i}-\bmu_{\alpha})(\Y_{\alpha i}-\bmu_{\alpha})^T\W(\overline{\Y}_{\beta(j)}-\bmu_{\beta})(\Y_{{\beta} j}-\bmu_{\beta})^T\W^{1/2}\bigg\},\\
&D_{4}=\frac{1}{n_\alpha n_\beta}\sum_{i=1}^{n_\alpha}\sum_{j=1}^{n_{\beta}}{\rm{tr}}\bigg\{\W^{1/2}(\Y_{\alpha i}-\bmu_{\alpha})\bmu_{\alpha}^T\W(\Y_{{\beta} j}-\bmu_{\beta})(\Y_{{\beta} j}-\bmu_{\beta})^T\W^{1/2}\bigg\},\\
&D_{5}=\frac{1}{n_\alpha n_\beta}\sum_{i=1}^{n_\alpha}\sum_{j=1}^{n_{\beta}}{\rm{tr}}\bigg\{\W^{1/2}(\Y_{\alpha i}-\bmu_{\alpha})(\Y_{\alpha i}-\bmu_{\alpha})^T\W(\Y_{{\beta} j}-\bmu_{\beta})\bmu_{\beta}^T\W^{1/2}\bigg\},\\
&D_{6}=\frac{-1}{n_\alpha n_\beta}\sum_{i=1}^{n_\alpha}\sum_{j=1}^{n_{\beta}}{\rm{tr}}\bigg\{\W^{1/2}(\overline{\Y}_{\alpha(i)}-\bmu_{\alpha})\bmu_{\alpha}^T\W(\Y_{{\beta} j}-\bmu_{\beta})(\Y_{{\beta} j}-\bmu_{\beta})^T\W^{1/2}\bigg\},\\
&D_{7}=\frac{-1}{n_\alpha n_\beta}\sum_{i=1}^{n_\alpha}\sum_{j=1}^{n_{\beta}}{\rm{tr}}\bigg\{\W^{1/2}(\Y_{\alpha i}-\bmu_{\alpha})(\Y_{\alpha i}-\bmu_{\alpha})^T\W(\overline{\Y}_{\beta(j)}-\bmu_{\beta})\bmu_{\beta}^T\W^{1/2}\bigg\},\\
&D_{8}=\frac{1}{n_\alpha n_\beta}\sum_{i=1}^{n_\alpha}\sum_{j=1}^{n_{\beta}}{\rm{tr}}\bigg\{\W^{1/2}(\overline{\Y}_{\alpha(i)}-\bmu_{\alpha})(\Y_{\alpha i}-\bmu_{\alpha})^T\\
&\quad\quad\W(\overline{\Y}_{\beta(j)}-\bmu_{\beta})(\Y_{{\beta} j}-\bmu_{\beta})^T\W^{1/2}\bigg\},\\
&D_{9}=\frac{-1}{n_\alpha n_\beta}\sum_{i=1}^{n_\alpha}\sum_{j=1}^{n_{\beta}}{\rm{tr}}\bigg\{\W^{1/2}(\Y_{\alpha i}-\bmu_{\alpha})\bmu_{\alpha}^T\W(\overline{\Y}_{\beta(j)}-\bmu_{\beta})(\Y_{{\beta} j}-\bmu_{\beta})^T\W^{1/2}\bigg\},\\
&D_{10}=\frac{-1}{n_\alpha n_\beta}\sum_{i=1}^{n_\alpha}\sum_{j=1}^{n_{\beta}}{\rm{tr}}\bigg\{\W^{1/2}(\overline{\Y}_{\alpha(i)}-\bmu_{\alpha})(\Y_{\alpha i}-\bmu_{\alpha})^T\W(\Y_{{\beta} j}-\bmu_{\beta})\bmu_{\beta}^T\W^{1/2}\bigg\},\\
&D_{11}=\frac{1}{n_\alpha n_\beta}\sum_{i=1}^{n_\alpha}\sum_{j=1}^{n_{\beta}}{\rm{tr}}\bigg\{\W^{1/2}(\overline{\Y}_{\alpha(i)}-\bmu_{\alpha})\bmu_{\alpha}^T\W(\overline{\Y}_{\beta(j)}-\bmu_{\beta})(\Y_{{\beta} j}-\bmu_{\beta})^T\W^{1/2}\bigg\},\\
&D_{12}=\frac{1}{n_\alpha n_\beta}\sum_{i=1}^{n_\alpha}\sum_{j=1}^{n_{\beta}}{\rm{tr}}\bigg\{\W^{1/2}(\overline{\Y}_{\alpha(i)}-\bmu_{\alpha})(\Y_{\alpha i}-\bmu_{\alpha})^T\W(\overline{\Y}_{\beta(j)}-\bmu_{\beta})\bmu_{\beta}^T\W^{1/2}\bigg\},\\
&D_{13}=\frac{1}{n_\alpha n_\beta}\sum_{i=1}^{n_\alpha}\sum_{j=1}^{n_{\beta}}{\rm{tr}}\bigg\{\W^{1/2}(\Y_{\alpha i}-\bmu_{\alpha})\bmu_{\alpha}^T\W(\Y_{{\beta} j}-\bmu_{\beta})\bmu_{\beta}^T\W^{1/2}\bigg\},\\
&D_{14}=\frac{-1}{n_\alpha n_\beta}\sum_{i=1}^{n_\alpha}\sum_{j=1}^{n_{\beta}}{\rm{tr}}\bigg\{\W^{1/2}(\overline{\Y}_{\alpha(i)}-\bmu_{\alpha})\bmu_{\alpha}^T\W(\Y_{{\beta} j}-\bmu_{\beta})\bmu_{\beta}^T\W^{1/2}\bigg\},\\
&D_{15}=\frac{-1}{n_\alpha n_\beta}\sum_{i=1}^{n_\alpha}\sum_{j=1}^{n_{\beta}}{\rm{tr}}\bigg\{\W^{1/2}(\Y_{\alpha i}-\bmu_{\alpha})\bmu_{\alpha}^T\W(\overline{\Y}_{\beta(j)}-\bmu_{\beta})\bmu_{\beta}^T\W^{1/2}\bigg\},\\
&D_{16}=\frac{1}{n_\alpha n_\beta}\sum_{i=1}^{n_\alpha}\sum_{j=1}^{n_{\beta}}{\rm{tr}}\bigg\{\W^{1/2}(\overline{\Y}_{\alpha(i)}-\bmu_{\alpha})\bmu_{\alpha}^T
\W(\overline{\Y}_{\beta(j)}-\bmu_{\beta})\bmu_{\beta}^T\W^{1/2}\bigg\}.
\end{align*}
From the above equation, we can conclude that $\E(D_{1})={\rm{tr}}(\W\bSigma_{\alpha}\W\bSigma_{\beta})$ and $\E(D_{l})=0$ for $l=2,\dots,16$.
Therefore, we have $\E\left({\rm{tr}}(\W\widehat{\bSigma_{\alpha}}\W\widehat{\bSigma_{\beta}})\right)={\rm{tr}}(\W\bSigma_{\alpha}\W\bSigma_{\beta})$, and now we will prove $\Var\left({\rm{tr}}(\W\widehat{\bSigma_{\alpha}}\W\widehat{\bSigma_{\beta}})\right)=o\left({\rm{tr}}^{2}(\W\bSigma_{\alpha}\W\bSigma_{\beta})\right)$.
Note that 
$$\Var\left({\rm{tr}}(\W\widehat{\bSigma_{\alpha}}\W\widehat{\bSigma_{\beta}})\right)\leq 16\sum_{l=1}^{16}\Var(D_{l}).$$
Since 
\begin{align*}
\Var(D_{1})&=\E(D_{1}^{2})-\E^{2}(D_{1})\\
&=\frac{1}{n_{\alpha}^{2} {n_\beta}^{2}}\E\bigg[\sum_{i=1}^{n_\alpha}\sum_{j=1}^{n_{\beta}}{\rm{tr}}\bigg\{\W^{1/2}(\Y_{\alpha i}-\bmu_{\alpha})(\Y_{\alpha i}-\bmu_{\alpha})^T\\
&\quad\W(\Y_{{\beta} j}-\bmu_{\beta})(\Y_{{\beta} j}-\bmu_{\beta})^T\W^{1/2}\bigg\}\\
&\quad\times \sum_{i_{1}=1}^{n_\alpha}\sum_{j_{1}=1}^{n_{\beta}}{\rm{tr}}\bigg\{\W^{1/2}(\Y_{\alpha i_{1}}-\bmu_{\alpha})(\Y_{\alpha i_{1}}-\bmu_{\alpha})^T\\
&\quad\W(\Y_{{\beta} j_{1}}-\bmu_{\beta})(\Y_{{\beta} j_{1}}-\bmu_{\beta})^T\W^{1/2}\bigg\}\bigg]-{\rm{tr}}^{2}(\W\bSigma_{\alpha}\W\bSigma_{\beta})\\
&=\frac{1}{n_\alpha n_\beta}\E\bigg[{\rm{tr}}\bigg\{\W^{1/2}(\Y_{\alpha i}-\bmu_{\alpha})(\Y_{\alpha i}-\bmu_{\alpha})^T\\
&\quad\W(\Y_{{\beta} j}-\bmu_{\beta})(\Y_{{\beta} j}-\bmu_{\beta})^T\W^{1/2}\bigg\}\bigg]^2\\
&\quad+\frac{(n_\beta-1)}{n_\alpha n_\beta}\E\bigg[(\Y_{\alpha i}-\bmu_{\alpha})^T\W(\Y_{{\beta} j}-\bmu_{\beta})(\Y_{{\beta} j}-\bmu_{\beta})^T\W(\Y_{\alpha i}-\bmu_{\alpha})\\
&\quad(\Y_{\alpha i}-\bmu_{\alpha})^T\W(\Y_{{\beta} j_{1}}-\bmu_{\beta})(\Y_{{\beta} j_{1}}-\bmu_{\beta})^T\W(\Y_{\alpha i}-\bmu_{\alpha})\bigg]\\
&\quad+o\left({\rm{tr}}^{2}(\W\bSigma_{\alpha}\W\bSigma_{\beta})\right).
\end{align*}
Through lemma \ref{lem 2.2}, lemma \ref{lem 2.3} and assumptions A1, A3, A4, we have 
\begin{align*}
&\frac{1}{n_\alpha n_\beta}\E\bigg[{\rm{tr}}\bigg\{\W^{1/2}(\Y_{\alpha i}-\bmu_{\alpha})(\Y_{\alpha i}-\bmu_{\alpha})^T\W(\Y_{{\beta} j}-\bmu_{\beta})(\Y_{{\beta} j}-\bmu_{\beta})^T\W^{1/2}\bigg\}\bigg]^2\\
&\leq \frac{1}{n_\alpha n_\beta}\bigg[(3+\Delta)\rm{tr}^{2}(\W\bSigma_{\alpha}\W\bSigma_{\beta})+(3+\Delta)(2+\Delta)\rm{tr}(\W\bSigma_{\alpha}\W\bSigma_{\beta})^{2}\bigg]\\
&=o\left({\rm{tr}}^{2}(\W\bSigma_{\alpha}\W\bSigma_{\beta})\right)
\end{align*}
and
\begin{align*}
&\frac{(n_\beta-1)}{n_\alpha n_\beta}\E\bigg[(\Y_{\alpha i}-\bmu_{\alpha})^T\W(\Y_{{\beta} j}-\bmu_{\beta})(\Y_{{\beta} j}-\bmu_{\beta})^T\W(\Y_{\alpha i}-\bmu_{\alpha})\\
&\quad(\Y_{\alpha i}-\bmu_{\alpha})^T\W(\Y_{{\beta} j_{1}}-\bmu_{\beta})(\Y_{{\beta} j_{1}}-\bmu_{\beta})^T\W(\Y_{\alpha i}-\bmu_{\alpha})\bigg]\\
&\leq \frac{(n_\beta-1)}{n_\alpha n_\beta}\bigg[\rm{tr}^{2}(\W\bSigma_{\alpha}\W\bSigma_{\beta})+(2+\Delta)\rm{tr}(\W\bSigma_{\alpha}\W\bSigma_{\beta})^{2}\bigg]\\
&=o\left({\rm{tr}}^{2}(\W\bSigma_{\alpha}\W\bSigma_{\beta})\right).
\end{align*}
Therefore, we have
$$\Var(D_{1})=o\left({\rm{tr}}^{2}(\W\bSigma_{\alpha}\W\bSigma_{\beta})\right).$$
By using a similar method, we can obtain 
$$\Var(D_{l})=o\left({\rm{tr}}^{2}(\W\bSigma_{\alpha}\W\bSigma_{\beta})\right),for\quad l=2,\dots,16.$$
Further, we get $\Var\left({\rm{tr}}(\W\widehat{\bSigma_{\alpha}}\W\widehat{\bSigma_{\beta}})\right)=o\left({\rm{tr}}^{2}(\W\bSigma_{\alpha}\W\bSigma_{\beta})\right).$ So 
$\frac{{\rm{tr}}(\W\widehat{\bSigma_{\alpha}}\W\widehat{\bSigma_{\beta}})}{{\rm{tr}}(\W\bSigma_{\alpha}\W\bSigma_{\beta})}\stackrel{p}{\longrightarrow}1.$
In summary, we can obtain $\frac{\widehat\sigma^{2}({\rm{T}}_{n_0})}{\sigma^{2}({\rm{T}}_{n_0})}\stackrel{p}{\longrightarrow}1.$ The proof of this theorem has been completed.

\subsection*{Proof of theorem 2.4}
\noindent Based on condition \eqref{eq2.26}, we have 
$$\frac{\S_n}{\sigma({\rm{T}}_{n_0})}=o_{p}(1).$$
Further, from Theorem \ref{th2.3} and \eqref{eq2.10}, we can derive 
\begin{align*}
P\left(\frac{{\rm{T}}_{n}}{\widehat\sigma({\rm{T}}_{n_0})}\geq z_{\alpha}\right)
&=P\left(\frac{{\rm{T}}_{n_0}+2\S_{n}+{\lVert\D_{\theta}\bmu\rVert}^2}{\sigma({\rm{T}}_{n_0})}\frac{\sigma({\rm{T}}_{n_0})}{\widehat\sigma({\rm{T}}_{n_0})}\geq z_{\alpha}\right)\\
&=\Phi\left(-z_{\alpha}+\frac{{\lVert\D_{\theta}\bmu\rVert}^2}{\sigma({\rm{T}}_{n_0})}\right)(1+o(1)) .
\end{align*}
Then the theorem proof is complete.

\end{appendices}

\end{document}